\documentclass[11pt]{article}   	                 
\usepackage{amsmath,amsthm, amssymb,amsfonts,amscd,mathabx}

\usepackage[colorlinks=true,linkcolor=blue,citecolor=red]{hyperref}
\usepackage[pdftex]{graphicx}                 
\usepackage{wrapfig}                               
\usepackage[font=footnotesize]{caption}                     
\usepackage{float}

\numberwithin{equation}{section}
\newtheorem{theorem}{Theorem}[section]

\newtheorem{proposition}[theorem]{Proposition}

\newtheorem{hypothesis}[theorem]{Hypothesis}
\newtheorem{remark}[theorem]{Remark}

\newenvironment{Proof}[1][.]%
{\begin{trivlist}\item[]\textbf{Proof#1 }}%
{\qed\end{trivlist}}

 {\begin{trivlist}\item[]\textbf{Acknowledgments }}{\end{trivlist}}

\newcommand{\R}{\mathbb{R}}
\newcommand{\C}{\mathbb{C}}

\newcommand{\Z}{\mathbb{Z}}

\newcommand{\Rmnum}[1]{\uppercase\expandafter{\romannumeral #1\relax}}

\newcommand{\rmO}{\mathrm{O}}

\newcommand{\rmd}{\mathrm{d}}
\newcommand{\rme}{\mathrm{e}}
\newcommand{\rmi}{\mathrm{i}}

\renewcommand{\leq}{\leqslant}
\renewcommand{\geq}{\geqslant}

\def\ker{\mathop{\mathrm{\,Ker}\,}}

\def\rg{\mathop{\mathrm{\,Rg}\,}}

\def\eps{\varepsilon}

\newcommand{\elll}{L}

\newfam\bifam
\font\tenbi=cmmib10 scaled \magstep1 \font\sevenbi=cmmib10 at 11pt
\font\fivebi=cmmib10 at 6pt \textfont\bifam = \tenbi
\scriptfont\bifam = \sevenbi \scriptscriptfont\bifam= \fivebi

\begin{document}

\begin{center}
{\fontsize{15}{15}\fontfamily{cmr}\fontseries{b}\selectfont{
Reversible switching due to attraction and repulsion: \\clusters, gaps, sorting, and mixing
}}\\[0.2in]
Arnd Scheel$\,^1$ and Angela Stevens$\,^2$\\[0.1in]
\textit{\footnotesize
$\,^1$ University of Minnesota, School of Mathematics, 206 Church St. S.E.,
Minneapolis, MN 55455, USA\\
$\,^2$University of M\"unster, Applied Mathematics, Einsteinstr. 62,
D-48149 M\"unster, Germany    
}
\end{center}

\begin{abstract}
\noindent
We describe a phase transition in continuum limits of interacting particle systems that exhibits a vertical bifurcation diagram. The transition is mediated by a competition short-range repulsion and long-range attraction. As a consequence of the transition, infinitesimal parameter variations allow switching between uniform distribution and clusters in single-species models, and between mixed and sorted states in multi-species contexts, without hysteresis. Our main technical contribution is a universal expansion for the size of vacuum bubbles that arise in this phase transition and a quantitative analysis of the effect of noise.
\end{abstract}
\section{Motivation and Setup}

We present here results on a curious reversible phase transition in continuum limits of interacting-particle systems. Motivation for this analysis stems from the apparent ability of biological systems to easily {but robustly} switch back and forth between different {functional} states through only small changes in the ambient environment. This switching is reversible, that is, the small change in the environment can be undone so that the system then reverts to the original state. This absence of hysteresis or memory in the system is unusual in more generic dynamical systems, where transitions between qualitatively different states often are represented as symmetry-breaking bifurcations.
\begin{figure}[b!]
\hfill\includegraphics[width=0.3\textwidth]{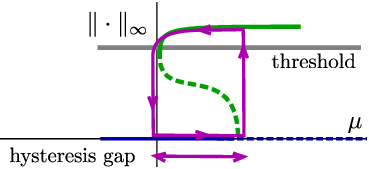}\hfill%
\includegraphics[width=0.3\textwidth]{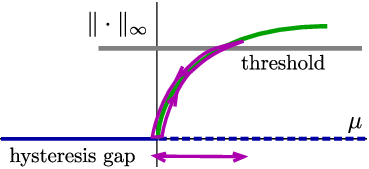}\hfill%
\includegraphics[width=0.3\textwidth]{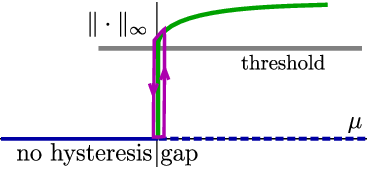}\hfill $ $%
\caption{Switching with finite hysteresis between states with difference of order one, and back, in the case of sub- (left) and supercritical pitchfork bifurcations (center), as opposed to switching in the interacting particle systems  (right) without hysteresis in the infinite-size limit.
}\label{f:1}
\end{figure}
Such bifurcations usually incur hysteresis upon switching between states that differ by a magnitude of order one, that is, parameter changes of order one are necessary to go through a cycle of changing the system to a new state, and back; see Fig. \ref{f:1}.
In this introduction, we first describe sorting and unsorting in lipid rafts,  a particular biological background that motivated this work, then introduce the general setup of short-range repulsion and medium-range attraction, and finally describe the main mathematical results.

\paragraph{Switching in lipid rafts.}

We are motivated by the concept of lipid-rafts in biological membranes
\cite{Simons1997}.
The { so-called}  ``raft hypothesis'' introduced there states that preferential interactions between
sterols and certain phospholipids can establish {cellular} membrane domains with
distinctive protein and lipid compositions  compared to the surrounding membrane.
So far, there is no conclusive evidence for the existence of rafts in vivo, which might be explained by the fact that coexisting ordered and disordered lipid phases in the membranes of living cells
may have specific, very small time and length scales and thus might be difficult to observe.
Studies of giant plasma membrane vesicles \cite{Baumgart} do support
the lipid-raft hypothesis though, at least to some extent. In these vesicles, a macroscopic phase separation
into ordered and disordered membrane parts can be observed; we refer to  \cite{LEVENTAL2020341} for a review. 
Our interest here is in how very basic mechanisms, namely attraction and repulsion among two species, can facilitate such observed transitions between order and disorder. In our model, we are particularly concerned with the robustness of the switching process: in a complex system with many unknown parameters and one specific control parameter, do we expect first to find a critical switching value, and do we then expect that switching is, as the name insinuates, reversible?

Our results do indeed point to such robust mechanisms that would enable fast and precise functional change in  biological membranes via
the formation of ordered and disordered raft nanodomains, and even gaps or channels, and their subsequent dissolution, all reversibly mediated by small manipulation of parameters in the ambient environment.

\begin{figure}
\hfill
\includegraphics[width=0.35\textwidth]{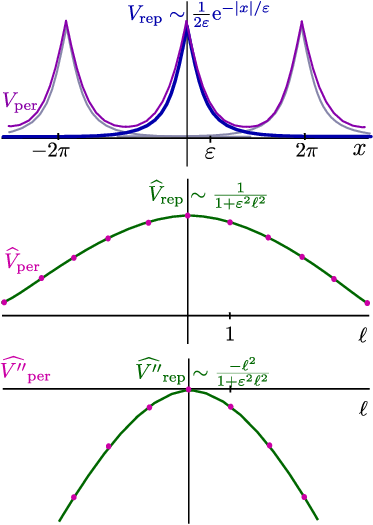}\hfill
\includegraphics[width=0.35\textwidth]{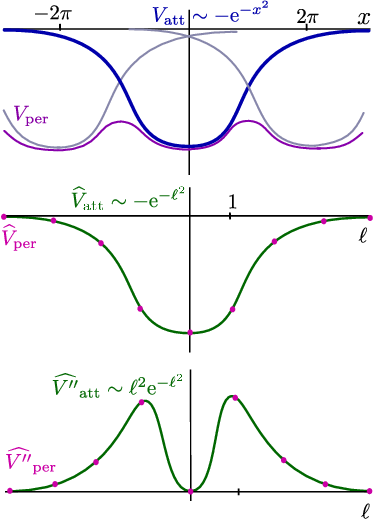}\hfill $ $
\caption{Schematics of repulsive (left) and attractive (right) potentials $V_\mathrm{rep}$ and $V_\mathrm{att}$ together with their periodized versions (top) obtained by adding all $2\pi \Z$-translates; Fourier transforms of potentials and their periodized version shown in the next row, where periodization corresponds to sampling; Fourier transform of the second derivative that determines stability is shown in the bottom row, illustrating how $\widehat{V}''_\mathrm{rep}+\mu \widehat{V}''_\mathrm{att}$ is positive precisely for $|\ell|=1$ past some critical value of $\mu\gtrsim \mu_*$, thus triggering the instability described in Hyp. \ref{h:ker}.}\label{f:pot}
\end{figure}

\paragraph{Setup: short-range repulsion,  long-range attraction, and the crystalline state.}
Our analysis starts  with the motion of $N$ interacting particles at one-dimensional positions $\underline{x}=(x_j)_{1\leq j\leq N}\in\R^N$. The interaction is described as an $\ell^2$-gradient flow to an interaction energy with pairwise interaction-potential $V$,
\begin{equation}\label{e:energy}
\mathcal{E}(\underline{x})=\sum_{1\leq j\leq N}\sum_{1\leq  m\leq N} V(x_j-x_m),
\end{equation}
that is,
\begin{equation}\label{e:particleode}
x'_j=-\nabla_{x_j}\mathcal{E}(\underline{x})=-\sum_{m\neq j}V'(x_j-x_m).
\end{equation}
To more naturally explain the setup, we assume that we can decompose $V$ into an attractive and a repulsive part $V=V_\mathrm{rep}+\mu V_\mathrm{att}$, both even and differentiable away from the origin,  so that
$V_\mathrm{rep}'(\xi)\xi<0$ and $V_\mathrm{att}'(\xi)\xi>0$. The parameter $\mu {> 0}$ measures the strength of attraction and will be our main bifurcation parameter. We think of the repulsive potential as strong but short-range, and the attractive potential as smooth but medium- or long-range; see Fig. \ref{f:pot}.
To fix ideas, assume that
\begin{equation}\label{e:gauss}
V_\mathrm{rep}(\xi)=\frac{1}{2\eta}\rme^{-|x|/\eta},\qquad \qquad V_\mathrm{att}(\xi)=-\rme^{-x^2},
\end{equation}
with $\eta\ll 1$. The precise form of these potentials is not all that relevant for our analysis except for the concavity of $V_\mathrm{rep}$ near the origin that precludes clustering at single points. In fact, one can easily see that for narrowly spaced particles, the concavity of $V_\mathrm{rep}$ precludes occupation of the same position by several particles: the energy of three equally spaced particles is strictly less than the energy of the same configuration with the middle particle moved towards either one of the outer particles by concavity.

A particular solution to the particle evolution \eqref{e:particleode} is the crystalline equilibrium, albeit with infinitely many particles $N=\infty$, where $x_j=\rho j$, $j\in\Z$, and the spacing $\rho>0$ is arbitrary. Linearizing at this state, we find
\begin{equation}\label{e:linparticle}
y_j'=-\sum_{m\neq j}
V''(\rho(j-m))(y_j-y_m),
\end{equation}
or, after Fourier transform $y_j=\int_0^{2\pi}\rme^{\rmi\sigma j}\hat{y}(\sigma)\rmd\sigma$,
\begin{equation}\label{e:lineparticlefourier}
\hat{y}'(\sigma)=\lambda(\sigma)\hat{y}(\sigma),\qquad
\lambda(\sigma)=\sum_{k\neq 0}V''(\rho k)(\rme^{-\rmi\sigma k}-1)= 
\frac{2\pi}{\rho}\sum_j \left(\widehat{V''}\left(\frac{2\pi j+\sigma}{\rho}\right)-\widehat{V''}\left(\frac{2\pi j}{\rho} \right)\right),
\end{equation}
where we used Poisson summation with the Fourier transform $\widehat{V''}(\ell)=\frac{1}{2\pi}\int V''(\xi)\rme^{-\rmi\ell\xi}\rmd\xi$. Note that we can allow for $V''$ to be singular at the origin since the summation does not include $k=0$.
For small $\rho$ and assuming that $\widehat{V}''(\ell)\to 0$ for $\ell\to\infty$, sufficiently rapidly, possibly after renormalizing by subtracting the singular part which does not appear in the summation, we see that the crystal is stable when $\widehat{V}''<0$.

\paragraph{The limit of infinite particle numbers.}
Letting $N\to\infty$ and rescaling time, one finds at least formally a continuum description for the density of particles $u(t,x)$, the Vlasov equation, which in our case reads
\begin{equation}\label{e:vlasov}
u_t=(u\cdot (V'*u))_x, \qquad (V*u)(x)=\int_y V(x-y)u(y)\rmd y;
\end{equation}
we refer to \cite{jabin} for a review, 


The crystalline state now corresponds to a state of uniform density, which we can normalize as $u(x)\equiv 1$. The linear stability of this uniform density is determined by the linearization
\begin{equation}\label{e:vlasovlin}
v_t=V''*v,\qquad  \text{ or, after Fourier transform, } \qquad\hat{v}_t= \widehat{V''}\cdot \hat{v}.
\end{equation}
The crystalline state then is stable precisely when the Fourier transform of $\widehat{V''}$ is non-positive, similar to our findings for the crystal.

\paragraph{Instabilities in periodic geometry.}
Assuming periodic particle arrangements, say $x_{j+M}=x_j+L$ for all $j$, one can effectively reduce the dynamics to $M$ ordinary differential equations for $x_j$, $0\leq j\leq M-1$, with periodized potentials $V_\mathrm{per}(\xi)=\sum_{j\in\Z}V(\xi+jL)$,
\begin{equation}\label{e:partper}
x'_j=-\sum_{m\neq j,0\leq m<M} V_\mathrm{per}'(x_j-x_m).
\end{equation}
Similarly, one finds the Vlasov equation \eqref{e:vlasov} and its linearization at the crystalline state \eqref{e:vlasovlin} with periodic boundary condition on the domain $x\in (0,L)$, with the periodized potential $V_\mathrm{p}$. Since the Fourier transform of $V_\mathrm{p}''$ is simply the Fourier transform of $V''$, sampled on the dual grid $\ell\in \frac{2\pi}{L}\Z$, we can readily determine stability.

In the example of short-range repulsion and medium-range attraction $V=V_\mathrm{att}+\mu V_\mathrm{rep}$, the Fourier transform of $V_\mathrm{att}$ decays slowly, $\widehat{V}_\mathrm{att}\sim \frac{1}{1+\eta^2 k^2}$, while the Fourier transform of the smooth attractive potential decays rapidly, in the example \eqref{e:gauss} as $\rme^{-k^2}$ so that we expect an instability when $\mu$ increases past a finite value $\mu_*$; see again Fig. \ref{f:pot}. In the following, we choose $L=2\pi$, absorbing the possibly different period into the potential $V$ by scaling.

\paragraph{Multi-species systems of particles.}
The above discussion can be readily repeated for systems with different particle types, or species, choosing positions $x_{p,j}$ for all particle types $1\leq p\leq P$, interaction potentials $V_{pq}$ between particles of type $p$ and type $q$, all split into a short-range repulsion and a smooth, potentially attracting long-range term. In the Vlasov-limit, we obtain again \eqref{e:vlasov}, where now $V'$ is a $P\times P$-matrix and $u\in \R^P$ denotes densities of the different species.  Multiplication is understood as  $(u\cdot (V'*u))_p=u_p\cdot (V'*u)_p$. Stability is determined by the eigenvalues of the matrix $\widehat{V''}(\xi)$: negative eigenvalues lead to instability. We will focus on an example of two species, with strong intraspecies repulsion $V_{\mathrm{rep},jj}$, $j=1,2$, and slightly weaker cross-species repulsion $V_{\mathrm{rep},12}$, allowing for particles of different species to cross each other; see \S\ref{s:le} for details.

\paragraph{The reversible phase transition.}
To formulate our setup, fix a base state and hence total average concentrations as  $u\equiv 1_P=(1,\ldots,1)\in\R^P$. We also focus on the setting with $2\pi$-periodic boundary conditions and drop the subscript referring to the periodization of the potential. We shall focus on the simplest case of a one-dimensional kernel of the linearization, dominated by the first component. Therefore, write $\mathcal{L}= \widehat{V}(\ell)\in\C^{P\times P}$ and decompose
\[
\mathcal{L}=
  \left(\begin{array}{cc}
    \mathcal{L}_{11}& \mathcal{L}_{1h}\\
    \mathcal{L}_{h1}&\mathcal{L}_{hh}
  \end{array}\right),
  \qquad \text{ acting on } \C\times \C^{P-1}.
\]

\begin{hypothesis}[Minimal kernel and crossing at bifurcation]\label{h:ker}
We assume that the Fourier transform $\mathcal{L}=\widehat{V}(\ell)\in\C^{P\times P}$ of $V$ has a 1d-kernel at Fourier wavenumber $|\ell|=1$ for $\mu=0$, with eigenvector $e_0$, normalized so that $|e_0|_\infty=1$, and trivial kernel for $|\ell|\neq 1$. Without loss of generality, we may assume that the maximal amplitude in the kernel is attained in the first component, that is, $e_{0,1}=1$, and  $|e_{0,j}|<1$ for $2\leq j\leq P$. We also assume a crossing condition,
\begin{equation}\label{e:crossing}
\frac{\rmd}{\rmd\mu}\left(\mathcal{L}_{11}^{-1}\mathcal{L}_{1h}\mathcal{L}_{hh}^{-1}\mathcal{L}_{h1}\right)=1, \text{ at } \mu=0.
\end{equation}
%
\end{hypothesis}
Clearly, reparameterizing $\mu$ we can satisfy \eqref{e:crossing} provided the left-hand side is nonzero.

\begin{hypothesis}[Regularity of interaction potentials]\label{h:smooth}
We assume that for the interaction potential $V$ we have
\[
\hat{V}(\ell)=V_\infty \ell^{-2\beta}(1+\rmO(\ell^{-2})), \qquad \text{for } \ell\to\infty,
\]
 for  $\beta\in[0,1]$ and some $V_\infty>0$.
\end{hypothesis}
For $\beta=0$, this implies that we can write $V(x)=V_\infty \delta(x)+V_1(x)$ with $V_1'$ being H\"older continuous. For $\beta=1$, we can write $V(x)=V_\infty G(x) +V_1(x)$, where $G(x)$ is the resolvent kernel of $\partial_{xx}$, that is,
\begin{equation}\label{e:G}
G(x)=\frac{\rme^{\frac{2 \pi
   -|x|}{\eta}}+\rme^{\frac{|x|}{\eta}}}{2 \left(\rme^{\frac{2
   \pi }{\eta}}-1\right) \eta},\qquad \eta=V_\infty^{-1/2},
\end{equation}
and $V_1'''(x)$ is H\"older continuous. Alternatively, in this latter case, $V$ possesses a $|x|$-singularity at $x=0$ and is otherwise smooth.

Our next main results considers equilibrium solutions to \eqref{e:vlasov},
\begin{equation}
u\cdot (V'_\mu*u)=0, \qquad u(x)=u(x+2\pi) \text{ for all } x.
\end{equation}
We will construct weak solutions to this equation, under the smoothness assumption on the kernel $V$, that are continuous and satisfy
\begin{equation}\label{e:weak}
u(x)=0,\qquad \text{ or }\qquad 
(V_\mu'*u)(x)=0.
\end{equation}
We also shall fix the average mass, that is, we require $\int_0^{2\pi}u= 2\pi 1_P$.

\begin{figure}
\includegraphics[height=2in]{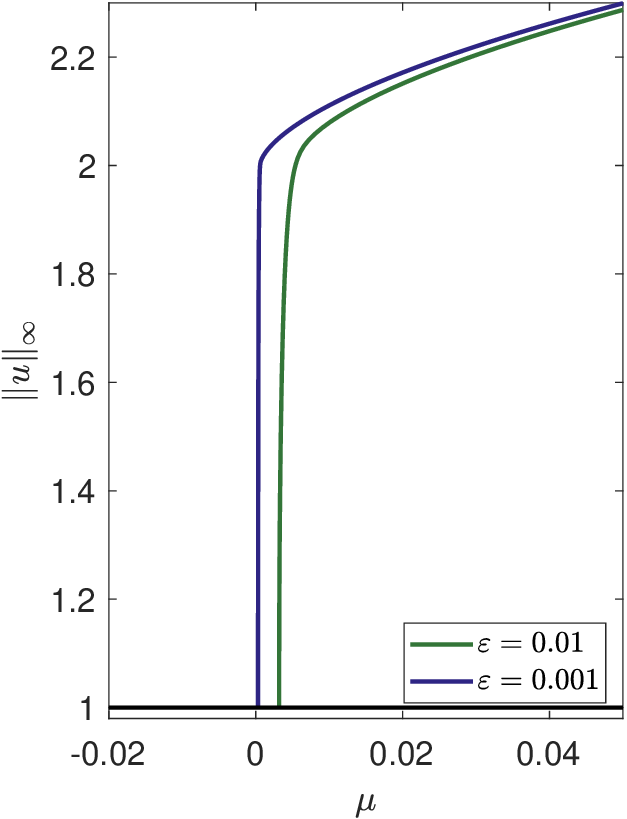}\hfill
\includegraphics[height=2in]{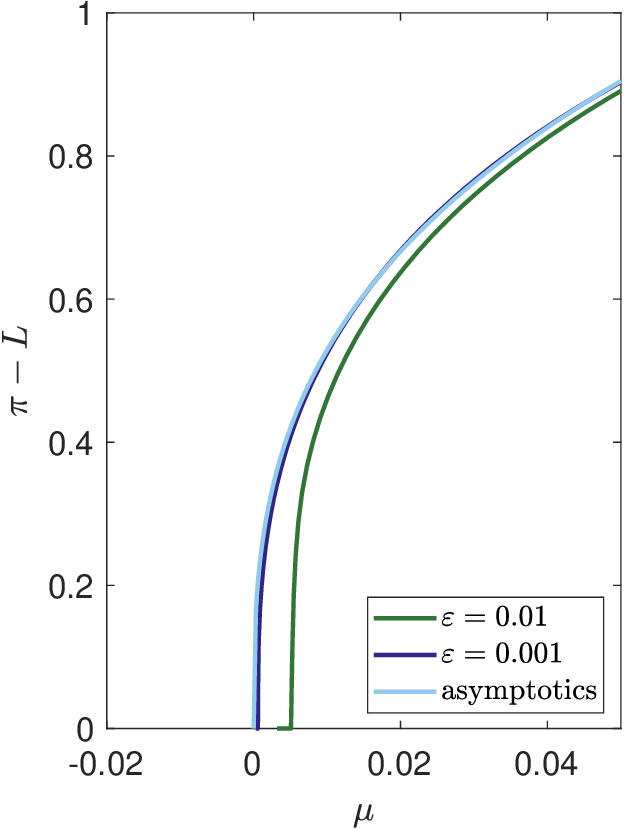}\hfill
\includegraphics[height=2in]{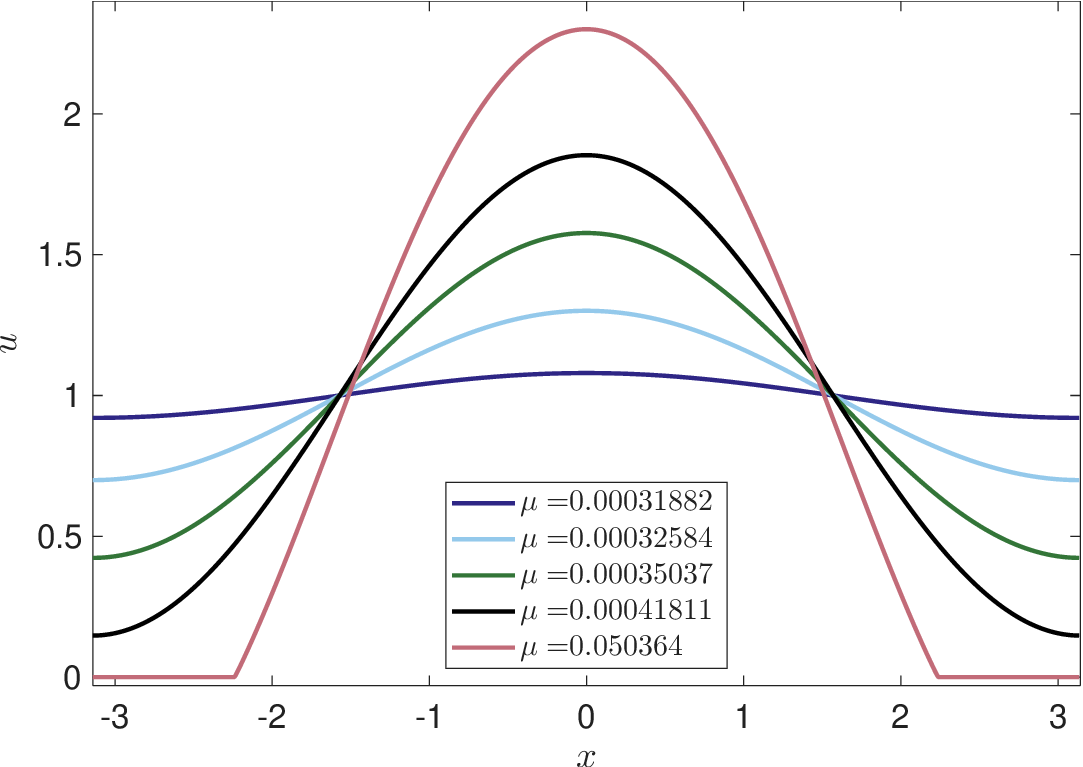}
\caption{Bifurcation of stationary solutions to \eqref{e:vlasov} with viscous regularization $\eps\partial_{xx}u$:  bifurcation diagram (left), bubble size expansion compared with asymptotics (center), and sample profiles ($\eps=0.001$, right). In all cases $V_\mathrm{att}(x)=(\mu-\frac{1}{\pi})\cos(x)+\frac{3}{20}\cos(2x)+\frac{1}{10} \cos(3x)$ and $V_\mathrm{rep}(x)=\delta(x)$.
}
\label{f:2a}
\end{figure}

\begin{theorem}[Vertical branch and gap formation]\label{t:1}
Assume Hypothesis \ref{h:ker} and \ref{h:smooth} for the interaction potential. We then have, for $0\leq\mu\ll 1$, a branch of even solutions of average mass $1_P$ consisting of
\begin{itemize}
\item a ``linear branch'' $u(x)=1_P+\rho e_0\cos(x)$, $0\leq \rho\leq 1$;
\item a ``small-vacuum branch'' $u(x;\mu)$, continuous in $0\leq \mu\ll 1$, uniformly in $x$, $u(x;0)=1_P+ e_0\cos(x)$, and $u_1(x;\mu)>0$ for $x$ in $(-\pi,\pi)$ precisely when $|x|<L(\mu)$ with $L(0)=\pi$, $L$ continuous; in the case of Dirac repulsion, we have the universal expansion
\begin{equation}\label{e:Luniv}
L(\mu)=\pi- \left(\frac{3\pi^2}{2}\right)^{1/3} \mu^{1/3}+\rmO(\mu^{2/3}).
\end{equation}
\end{itemize}
All other components $u_j(x;\mu)$, $2\leq j\leq P$ are strictly positive, by continuity.
\end{theorem}
Sample bifurcation diagrams that illustrate the result are shown in Fig. \ref{f:2a}. We emphasize that the exponent $1/3$ enhances the effect of switching: even very small changes of $\mu$ open vacuum regions of considerable size $\sim \mu^{1/3}$.

\begin{remark}[Universality of expansion]
The expansion is universal in the sense that the exponent $1/3$ and the precoefficient do not depend on the attractive potential. Of course, reparameterization of $\mu$ would change the precoefficient, so that the precise form here depends on the speed of crossing, normalized to 1 in \eqref{e:crossing}. One may wish to {compare} this to say a pitchfork bifurcation where the amplitude scales with $\sqrt{\mu}$ and a precoefficient that depends on the speed of crossing of the eigenvalue and, crucially, the coefficient of a cubic term, sometimes referred to as a Landau coefficient. In our case, the expansion does not depend on such a nonlinear term.
We do note however that we noticed a strong dependence of the gap size on the regularity of the repulsive potential; see Fig. \ref{f:frac}, below. It would be interesting to describe more generally the power law asymptotics of bubble sizes for different regularity classes of $V_\mathrm{rep}$.
\end{remark}

\paragraph{Comparison with prior work.}
There is an extensive literature on competing attraction and repulsion in interacting particle systems, although little from the perspective of dynamics and bifurcation taken here. The early work \cite{TBL} identifies a transition to patterned states with vacua as found here. The specific bifurcation diagram we show here was noted in the particular situation of a Dirac-$\delta$ repulsive potential and a Bessel-repulsive potential, motivated by one-dimensional chemotaxis in
\cite{berthelin} and again in \cite{ccww}. Existence of a linear branch at the transition in the qualitative shape of energy minimizers much more generally was noted in \cite{cky}. The systems case was, at least for two species, analyzed for instance in \cite{bdfs18,bcps20}, establishing existence of clustered and segregated minimizers, albeit not from the bifurcation perspective taken here.
Our results demonstrate that the observations for the specific chemotaxis setting are much more generally valid, for systems and for weaker repulsive potentials, at least locally near the bifurcation point, where we also give a more detailed description of the vacuum compared to \cite{berthelin,ccww}.

Turning back to the question of spatial organization on lipid rafts, we also mention here  \cite{GKRR}, where a mathematical model 
that couples bulk-diffusion to a Cahn-Hilliard type equation on the membrane. The latter model is well established as a {thermodynamically} valid description of phase transitions and possesses a macroscopically, ``built-in'' tendency to separate phases. We emphasize however that phase transitions modeled in this fashion are typically not reversible, with generic sub- and supercritical bifurcations that exhibit hysteresis.
%


\paragraph{Acknowledgments.} A.Sc. gratefully acknowledges support through grant NSF DMS-2205663 and generous hospitality at the University of M\"unster where part of this work was carried out. 
A. St. was supported by the DFG
(German Research Foundation) under Germany's Excellence Strategy
EXC 2044-390685587, 
Mathematics M\"unster: 
Dynamics - Geometry - Structure
and by an Ordway Visiting Professorship during a scientifically inspiring one month research stay in 
the welcoming atmosphere of the University
of Minneapolis.

\section{The linear branch and examples}\label{s:le}
Solving for equilibria of \eqref{e:vlasov} in a periodic geometry $x\in [0,2\pi)$, we integrate the right-hand side, find that the constant of integration is necessarily zero, to obtain
\begin{equation}\label{e:eq}
u\cdot (V'*u)=0.
\end{equation}
For equilibria without vacuum, that is, $u(x)>0$ for all $x$, this is equivalent to 
\begin{equation}\label{e:eq1}
\int_0^{2\pi}V(x-y)u(y)=\rho,
\end{equation}
for some $\rho\in\R^P$. With $V(\xi)=\sum_{\ell\in\Z} \widehat{V}_\ell \rme^{\rmi\ell x}$, a basis of solutions to this equation is given by functions $\rme^{\rmi\ell x}e_{\ell,j}$, where $\ell$ is such that $\widehat{V}_\ell$ is not invertible and the $e_{\ell,j}$ span its kernel, up to the trivial constant elements. Our assumptions then guarantees that at the bifurcation point, $\mu=0$, there is a one-dimensional space of even solutions with average density $1_P$, $u(x)=1_P+\rho\cos(x)e_0$, $\rho\in\R$; see Hypothesis \ref{h:ker}. Clearly, requiring $u(x)\geq 0$ restricts $|\rho|\leq 1$, recalling that we normalized $|e_0|_\infty=1$. Note also that changing the sign of $\rho$ {amounts simply} to shifting the solution by $\pi$ and we therefore have, up to translations, a unique interval of solutions $\rho\in[0,1]$. The analysis { in} the next sections explores a vicinity of the solution at $\rho=1$, $u_*(x)=1_P+\cos(x) e_0$.

\paragraph{Examples.}
The simplest examples arise for scalar equations, $P=1$, with for instance Dirac-$\delta$ repulsive potential and a cosine attractive potential, 
\[
u_t=(u(u +  V_\mu*u)_x)_x, \qquad V_\mu(x)=-\left(\frac{1}{\pi}+\mu\right)\cos(x).
\]
Here, the cosine is just the simplest periodic attractive potential and the model can be thought of as neglecting harmonics. We shall see below that in this case solutions with vacuum are given explicitly as the positive part of the sum of constants and a multiple of $\cos(x)$. Smoothed versions of the repulsive potential lead to the generalization
\[
u_t=(u((1-\eta^2\partial_{xx})^{-\beta}u +  V_\mu*u)_x)_x, \qquad V_\mu(x)=-(\mu_*+\mu)\cos(x),\quad 
\mu_*=\frac{1}{\pi(1+\eta^2)^\beta}.
\]
For $P=2$, with species $u_1$ and $u_2$, we focus on Dirac-repulsion and cosine attractive potentials,
\begin{equation}
(V_\mathrm{rep}*u)(z)=
  \left(
    \begin{array}{cc}
      a_{11}& a_{12}\\
      a_{21}& a_{22}
    \end{array}
  \right)
  \left(
    \begin{array}{c}
      u_{1}(z)\\
      u_{2}(z)
    \end{array}
  \right),
  \qquad 
  (V_\mathrm{att}*u)(z)=-
  \left(
    \begin{array}{cc}
      b_{11}& b_{12}\\
      b_{21}& b_{22}
    \end{array}
  \right)
  \left(
    \begin{array}{c}
      c_{1}(z)\\
      c_{2}(z)
    \end{array}
  \right),
\end{equation}
where $c_j(z)=\int_\xi \cos(z-\xi)u_j(\xi)\rmd\xi$, resulting in 
\begin{equation}\label{e:sys}
\begin{aligned}
u_{1,t}&=(u_1(a_{11}u_1+a_{12}u_2+b_{11} c_1+b_{12} c_2)_x)_x,\\
u_{2,t}&=(u_2(a_{21}u_1+a_{22}u_2+b_{21} c_1+b_{22} c_2)_x)_x.
\end{aligned}
\end{equation}
We assume repulsion $a_{ij}>0$ and strict ellipticity for positive masses, that is
\[
a_{11}a_{22}>a_{12}a_{21}.
\]
Inter-species repulsion then is weaker than intra-species repulsion, on average. The attractive potential is actually attractive if $b_{jm}>0$ for all $j,m$.  Linearizing at constant mass $1$ in all species, when $u=1_2$, and using Fourier transform, we find stability for all Fourier modes $|\ell|\neq 1$, and criticality at $\ell=1$ when
\[\mathrm{det}\,\left(
    \begin{array}{cc}
      a_{11}-\pi b_{11}& a_{12}-\pi b_{12}\\
      a_{21}-\pi b_{21}& a_{22}-\pi b_{22}
    \end{array}
 \right)=0. 
\]
To interpret two relevant cases, assume that $a_{11}-\pi b_{11}>0$, and $a_{22}-\pi b_{22}>0$, {that} is, ignoring effects of the other species, a species alone would not cluster. There are then two cases,
\begin{itemize}
\item[(JC)] \emph{joint clustering:} $a_{12}-\pi b_{12}<0$,
\item[(S)] \emph{segregation:} $a_{12}-\pi b_{12}>0$.
\end{itemize}
For joint clustering, the eigenvector $e_0=(e_{0,1},e_{0,2})^T$ in the kernel has $e_{0,1}\cdot e_{0,2}>0$, that is, concentrations in the kernel $u_1=1+\rho e_{0,1}\cos x$ and $u_2=1+\rho e_{0,2}\cos x$ exhibit maxima at the same {location.} In the case of segregation, $e_{0,1}\cdot e_{0,2}<0$ and maxima of concentrations of $u_1$ and $u_2$ are located at maximal distance $\pi$. This is of course intuitively plausible since $a_{12}-\pi b_{12}<0$ indicates inter-species effective attraction, beating out intra-species effective repulsion, whereas the other sign indicates strong inter-species repulsion overcoming the intra-species repulsion.

\begin{remark}[Vacuum formation]
The relative size of $|e_{0,1}|$ and $|e_{0,2}|$ determines which species forms a vacuum near criticality. Clearly, this relative size also depends on the size of the constant concentrations $u_1$ and $u_2$ which we normalized both to $1$. Our assumptions guarantee that these magnitudes differ,  $|e_{0,1}|\neq|e_{0,2}|$, so that we only analyze vacuum regions in one species.
\end{remark}

\begin{remark}[Sorting in one space dimension]
Inspecting the paths of particles in the discrete case $N<\infty$, the ordering of particles is preserved in scalar equations, $P=1$. Clearly, sorting and segregation for systems cannot preserve order and necessitates switching of position between different species. How such change of positions is effectuated in higher-dimensional systems is studied in some simple situations probabilistically, see for instance \cite{quastel}. If and how a one-dimensional model or a higher-dimensional generalization as the one studied here can capture this process does not seem to be understood.
\end{remark}

\section{Vacuum bubbles --- scalar rank-one potentials}
We wish to solve for equilibria of  \eqref{e:vlasov}  allowing for regions where $u$ vanishes. We therefore start with a very simple scalar case, where 
\[
V_\mathrm{rep}(x)=\delta(x),\qquad V_\mathrm{att}(x)=-\left(\frac{1}{\pi}+\mu\right)\cos(x).
\]
The equation for even equilibria with support in $-[\elll,\elll]$ then becomes
\begin{equation}\label{e:vac1}
\begin{aligned}
 u(x)-\left(\frac{1}{\pi}+\mu\right)\cos(x)\int_{-\elll}^\elll \cos(y)u(y)\rmd y-\rho&=0,\qquad |x|\leq \elll, \\
 \frac{1}{2\pi}\int_{-\elll}^\elll u(y)\rmd y-1&=0, \\
 u(\elll)&=0.
\end{aligned}
\end{equation}
for some $\rho\in \R$.
This integral operator is of rank 1 and we find from the first equation that $u(x)$ is necessarily given explicitly as $(A_0+A_1\cos(x))_+=\max(0,A_0+A_1\cos(x))$. This form of $u$ implicitly defines the boundary of the vacuum region as
$
 \elll=\arccos(-A_0/A_1),
$
thus solving the third equation in \eqref{e:vac1}.

Resolving the square root singularity of $\arccos(-A_0/A_1)$ near $A_0/A_1=1$, we scale
\[
 A_0=1+a_3 z_1^3,\qquad A_1=1+z_1^2+z_1^3,
\]
for $z_1$ small and $a_3=\rmO(1)$, to be determined as functions of the parameter $\mu$.

\begin{figure}
\hfill
\includegraphics[height=2in]{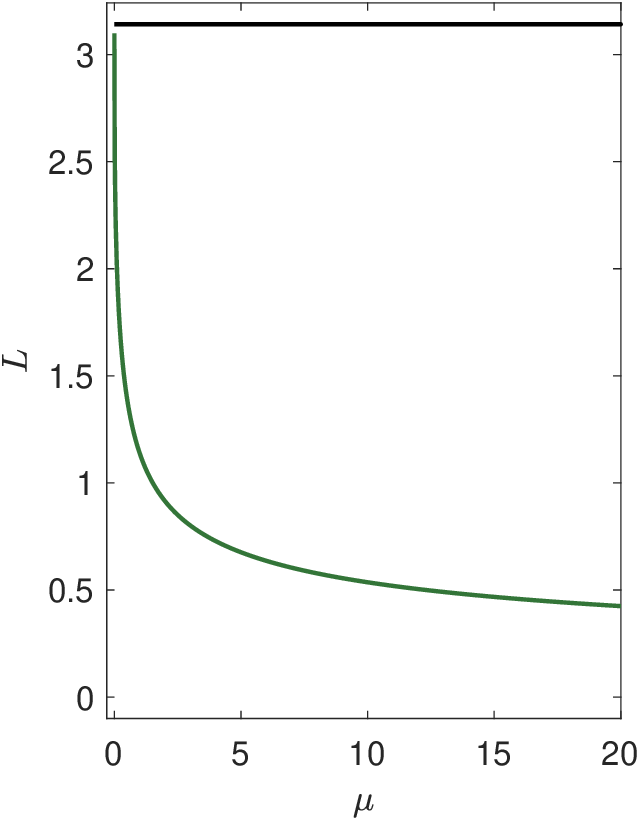}\hfill
\includegraphics[height=2in]{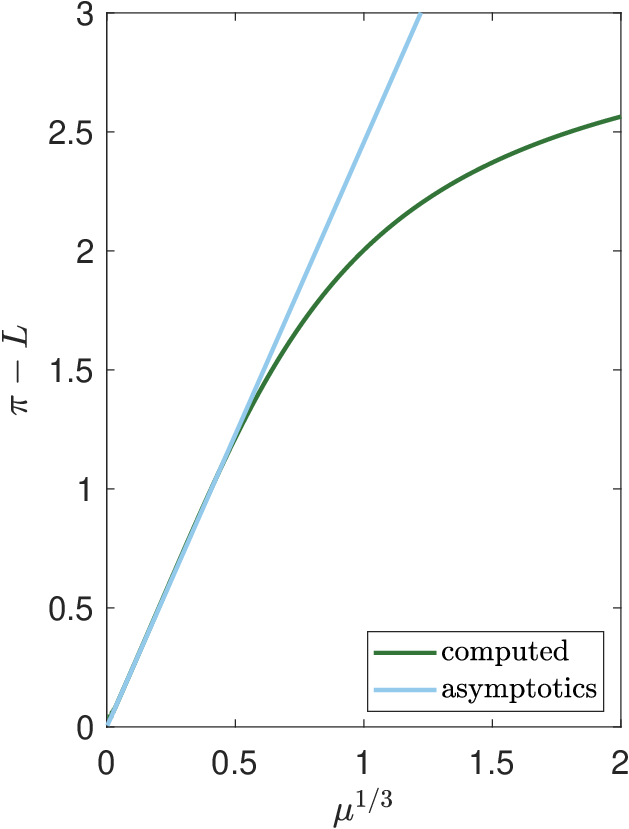}\hfill
\includegraphics[height=2in]{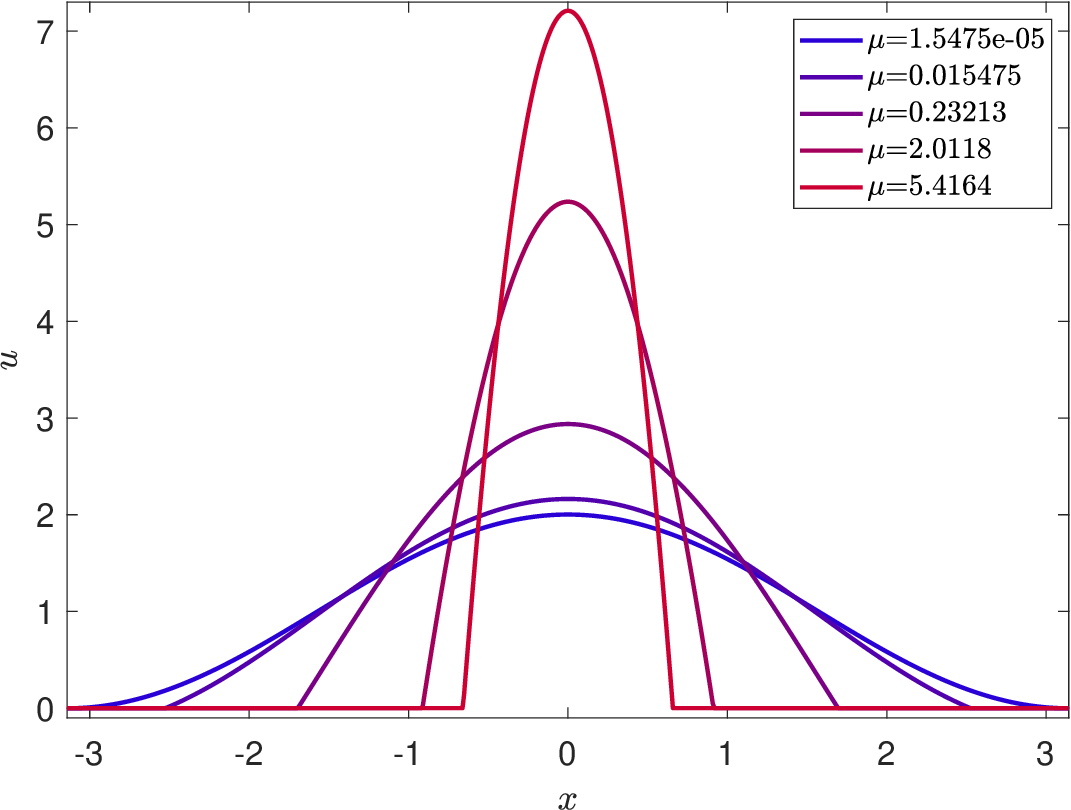}\hfill $ $
\caption{Global bifurcation diagram in the rank-one case $V(x)=\delta(x)-\mu\cos(x)$, computed by numerically solving for $(A_0,A_1)$ at fixed $\mu$. Size of the support $(-L,L)$ versus $\mu$ (left), comparison with asymptotic theory for small $\mu$ (center) and sample profiles (right).}
\label{f:lmu}
\end{figure}

The first equation consists of terms constant in $x$ and a multiple of $\cos(x)$. 
We then evaluate the constant terms in the first equation of \eqref{e:vac1} to find $\rho$ in terms of the remaining variables. The coefficient of $\cos(x)$ in the first equation of \eqref{e:vac1}, and the second equation in \eqref{e:vac1} evaluate to
\begin{align}
 \left(\frac{4\sqrt{2}}{3}+2a_3\pi\right) z_1^3 + \rmO(z_1^4)&=0,\label{e:coscomp}\\
 \mu\pi - \frac{4\sqrt{2}}{3\pi}z_1^3+\rmO(|\mu z_1^3|+z_1^4)&=0,\label{e:masscomp}
\end{align}
respectively. Dividing \eqref{e:coscomp} by $z_1^3$ and solving with the implicit function theorem for $a_3$, and then solving \eqref{e:masscomp} for $\mu$ with the implicit function theorem, we find the expansion
\begin{equation}\label{e:zsol}
 a_3=-\frac{4\sqrt{2}}{6\pi}+\rmO(z_1),\qquad \mu=\frac{4\sqrt{2}}{\pi^2}z_1^3+\rmO(z_1^4).
\end{equation}
In particular, using that the boundary of the vacuum region is given by $\elll=\sqrt{2}z_1+\rmO(z_1^3)$, and solving \eqref{e:masscomp} for $z_1$ as a function of $\mu$, we obtain
\begin{equation}\label{e:vacexp}
 \elll=\pi-\left(\frac{3\pi^2}{2}\right)^{1/3} \mu^{1/3}+\rmO(\mu^{2/3}).
\end{equation}
\begin{remark}[Large $\mu$]\label{r:lmu}
Trying to solve the effectively 4 equations \eqref{e:vac1}  $(A_0, A_1,L,\rho)$ for finite $\mu$, non-perturbatively, is difficult due to the presence of the transcendental term $\arccos(-A_0/A_1)$. Nevertheless, they can be solved numerically in a straightforward fashion. Results are illustrated in Fig. \ref{f:lmu}. The graph shows in particular that the asymptotics yield good predictions until the size of the vacuum is roughly half the size of the domain. It would be interesting to establish qualitative properties of this diagram analytically. One finds numerically that at leading order $A_0\sim -A_1 \sim \mu$ and $L\sim \mu^{-1/3}$ as $\mu\to\infty$. We emphasize that the expansion and shape of solutions for small $\mu$ is universal, roughly independent of the choice of the potential, so that the rank-one approximation from this section yields excellent predictions. For large $\mu$, however, this is of course not true and harmonics $\cos(m x)$, $m>1$, are relevant, leading potentially to secondary bifurcations.
\end{remark}
\begin{remark}[Finite rank]\label{r:fr}
The procedure here can in principle be adapted to solve for kernels $V$ given as finite Fourier polynomials since in this case the integral operator again has a finite-dimensional range. We may therefore choose $u$ inside this range parameterized by amplitudes of harmonics $A_j\cos(jx)$, before taking the positive part. This gives a finite-dimensional system of equations which one can again solve perturbatively or numerically. In the next section, we outline a somewhat more robust and functionally analytic sound approach for cases when $V$ is not simply a Fourier polynomial.
\end{remark}

\section{Vacuum bubbles --- general scalar potentials and a free-boundary formulation}\label{s:ift}

We first introduce a free-boundary formulation to solve for equilibria in the case where $V_\mathrm{rep}=\delta$ is the Dirac-$\delta$, and then show modifications for more general repulsive kernels. 

\subsection{Free boundary formulation of bubble-formation}\label{s:fbdy}

We rescale space to $x=\elll z/\pi$, write $v(z)=u(\elll z/\pi)$, and introduce the rescaled potential 
\[
 V_\mu^\elll(z)=\frac{\elll}{\pi}V_\mu(\frac{\elll}{\pi} z).
\]
The equation for equilibria in \eqref{e:vlasov}  then reads
\begin{align}
 0=F_v(v,\elll,\rho,\mu)[z]&:=v(z)-\int_{-\pi}^\pi V_\mu^\elll(z-\xi)v(\xi)\rmd\xi -\rho,\label{e:Fv}\\
 0=F_\mathrm{bc}(v)&:=v(\pi),\label{e:Fbc}\\
 0=F_m(v,\elll)&:=\frac{\elll}{2\pi^2}\int_{-\pi}^\pi v(\xi)\rmd\xi -1.\label{e:Fm}
\end{align}
We consider this as an equation with variables $(v,\elll,\rho,\mu)\in C^0([-\pi,\pi])\times \R^3$,
\[
 F=(F_v,F_\mathrm{bc},F_m):C^0([-\pi,\pi])\times \R^3\to C^0([-\pi,\pi])\times \R^2.
\]
One easily checks that the equation is smooth in all variables given sufficient smoothness in $V$. 

To further refine the analysis, we also introduce the projections
\begin{equation}\label{e:proj}
 P_0v=\frac{1}{2\pi}\int_{-\pi}^\pi v(\xi)\rmd\xi,\qquad P_1v=\frac{1}{\pi}\int_{-\pi}^\pi \cos(\xi)v(\xi)\rmd\xi\cos(z),\qquad P_\mathrm{h}v=\mathrm{id}-P_0-P_1,
\end{equation}
the associated functionals $\bar{P}_0v=\frac{1}{2\pi}\int_{-\pi}^\pi v(\xi)\rmd\xi$ and $\bar{P}_1v=\frac{1}{\pi}\int_{-\pi}^\pi \cos(\xi)v(\xi)\rmd\xi$, and decompose our solution to
\begin{equation}\label{e:vdecomp}
 v(z)=A_0+A_1\cos z + v_\mathrm{h},\qquad A_0=\bar{P}_0v,\quad A_1=\bar{P}_1v, \quad v_\mathrm{h}=P_\mathrm{h}v.
\end{equation}
{Substituting} this form of $v$ into \eqref{e:Fv}--\eqref{e:Fm} and decomposing \eqref{e:Fv} with the projections $P_{0,1,\mathrm{h}}$, we obtain a system
\begin{align}
F_v^0(A_0,A_1,v_\mathrm{h},\elll,\rho,\mu)&=0\nonumber\\
F_v^1(A_0,A_1,v_\mathrm{h},\elll,\rho,\mu)&=0\nonumber\\
F_v^\mathrm{h}(A_0,A_1,v_\mathrm{h},\elll,\rho,\mu)&=0\nonumber\\
F_\mathrm{bc}(A_0,A_1,v_\mathrm{h})&=0 \nonumber\\
F_m(A_0,\elll)&=0 \label{e:Fdec}
\end{align}
where we slightly abused notation for the new functions  $F_{\mathrm{bc}}$ and $F_\mathrm{m}$. We find explicitly
\begin{align}
F_\mathrm{bc}(A_0,A_1,v_\mathrm{h})&=A_0-A_1+v_\mathrm{h},\nonumber\\
F_m(A_0,\elll)&=\frac{\elll}{\pi}A_0-1,\label{e:bcm}
\end{align}
and,
\begin{align}
 F_v^1(A_0,A_1,v_\mathrm{h},\elll,\rho,\mu)&=
 A_1  - \frac{1}{\pi}\int_{-\pi}^\pi
\int_{-\pi}^\pi \cos(z)V_\mu^\elll(z-\xi)
(A_0+A_1\cos\xi + v_\mathrm{h}(\xi))
\rmd\xi\rmd z,
 \nonumber\\
 F_v^\mathrm{h}&=F_v-F_v^1 \cos z -F_v^0.\label{e:Fv01h}
\end{align}
We have the trivial solution, at the edge of the vertical branch, 
\[
 A_0^*=1,\qquad A_1^*=1,\qquad v_h^*=0,\qquad \elll_*=\pi,\qquad \rho_*=1,\qquad \mu_*=0.
\]
Writing $X_\mathrm{h}=\rg P_\mathrm{h}$, and reordering variables and components of the nonlinear function, we have a map
\begin{align}
 F=(F_v^h,F_v^1,F_m,F_\mathrm{bc}):X_\mathrm{h}\times \R^5&\to X_\mathrm{h}\times \R^4\nonumber\\
 (v_\mathrm{h},\rho,\mu,\elll,A_1,A_0)&\to (F_v^h,F_v^1,F_m,F_\mathrm{bc}),\label{e:Fcomp}
\end{align}
and $F(0,1,0,\pi,1,1)=0$. The linearization can be written as a $5\times6$-matrix, with first row and columns assumed to have range and domain in $X_\mathrm{h}$, respectively:
\begin{equation}\label{e:matrix}
 DF(0,1,0,\pi,1,1)=
 \left(\begin{array}{cccccc}
 \mathcal{L}_\mathrm{h} & 0& 0& g(z) & 0 & 0\\
 0& -1&  0 &\frac{1}{\pi} & 0 & 1\\
 0& 0& -\pi& 0 & 0 & 0 \\
 0 & 0 & 0 & \frac{1}{\pi}& 0 & 1 \\
 \delta(\cdot-\pi) & 0 & 0 & 0 & -1 & 1
 \end{array}\right),
\end{equation}
where $\delta(\cdot)$ is the Dirac-$\delta$ in the dual of $X_\mathrm{h}$,
\[
 [\mathcal{L}_\mathrm{h}v](z)=v(z)-[P_\mathrm{h}\int_{-\pi}^\pi V_\mu^\elll(\cdot-\xi)v(\xi)\rmd\xi](z)=v(z)-\int_{-\pi}^\pi \sum_{k=2}^\infty v_k\cos(k(z-\xi))v(\xi)\rmd\xi,
\]
and 
\begin{equation}\label{e:g}
 g(z)=-\int_{-\pi}^\pi \partial_\elll|_{\elll=\pi} V_0^\elll(z-\xi)(1+\cos(\xi))\rmd\xi.
\end{equation}
In the specific case of $V_0(z)=\frac{1}{\pi}\cos(z)$, we find $\mathcal{L}_\mathrm{h}=\mathrm{id}$, and $g(z)=\frac{1}{2\pi}(\cos z + 2z\sin z - 2)$, $g(\pi)=-\frac{3}{2\pi}$.

By the assumption that $v_k<\frac{1}{\pi}$, we have that $\mathcal{L}_h$ has a trivial kernel and thereby, as a compact perturbation of the identity, is bounded invertible.

Next, we wish to mimic row elimination to transform \eqref{e:matrix} into an upper triangular matrix. We therefore define the functional
\[
 b^*:=(\mathcal{L}^*)^{-1}\delta(\cdot-\pi)\in (X^\mathrm{h})^*
\]
and confirm that $b^*\mathcal{L}_\mathrm{h}=\delta(\cdot-\pi)$. We then redefine the last component of $F$ in \eqref{e:Fcomp} through
\begin{equation}\label{e:newF}
 \tilde{F}_\mathrm{bc}:=F_\mathrm{bc}-\langle b^*,F_v^\mathrm{h}\rangle +\pi \langle b^* ,g\rangle F_m,
\end{equation}
which yields the new linearization in upper triangular form,
\begin{equation}\label{e:matrix2}
 D\tilde{F}(0,1,0,\pi,1,1)=
 \left(\begin{array}{cccccc}
 \mathcal{L}_\mathrm{h} & 0& 0& g(z) & 0 & 0\\
 0& -1&  0 &\frac{1}{\pi} & 0 & 1\\
 0& 0& -\pi& 0 & 0 & 0 \\
 0 & 0 & 0 & \frac{1}{\pi}& 0 & 1 \\
 0 & 0 & 0 & 0 & -1& 1+\langle b^*,g\rangle
 \end{array}\right),
\end{equation}
with invertible diagonal entries. We conclude that $D_{(v_\mathrm{h},\rho,\mu,\elll,A_1)}\tilde{F}$ is invertible at  $(0,1,0,\pi,1,1)$ and we may therefore apply the implicit function theorem to find a unique solution
\begin{equation}\label{e:iftsol}
 (v_\mathrm{h},\rho,\mu,\elll,A_1)=(v_\mathrm{h},\rho,\mu,\elll,A_1)(A_0),
\end{equation}
for $|A_0-1|\ll 1$ sufficiently small. Note that, of course, the solution is only relevant when $\elll<\pi$. We therefore notice that, from $F_m$, $\elll=\pi/ A_0$, so that our solution is meaningful for $A_0\gtrsim 1$.

The expansion of the solution can now be found as usual with the implicit function theorem, although calculations are somewhat elaborate. We give some details in the following.

We define $\nu$ through  $A_0=1+\nu$, ${\nu}\gtrsim 0$, and expand formally
\[
  (v_\mathrm{h},\rho,\mu,\elll,A_1)=(0,1,0,\pi,1) + \sum_{j=1}^{3} (v_{\mathrm{h},j},\rho_j,\mu_j,\elll_j,A_{1,j})\nu^j + \rmO(\nu^4),
\]
substitute into $F=0$, and expand in $\nu$. At zeroth order we confirm that the base solution solves the equation. At first order, we recover the linear equation
$D\tilde{F}(0,1,0,\pi,1,1) (v_{\mathrm{h},1},\rho_1,\mu_1,\elll_1,A_{1,1},1)=0$, with solution
\[
 A_{1,1}=-\frac{1}{2},\ \elll_1=-\pi,\ \mu_1=0,\ \rho_1=\frac{1}{2}, \ v_{\mathrm{h},1}(\xi)=-1+\frac{1}{2}\cos\xi+\xi\sin\xi.
\]
At second order, we find a linear equation in $(v_{\mathrm{h},2},\rho_2,\mu_2,\elll_2,A_{1,2})$ with right-hand side given by evaluating nonlinearities at order $\nu^2$. Solving, we find
\[
 A_{1,2}=\frac{1}{12}(3+4\pi^2),\ \elll_2=\pi,\ \mu_2=0, \ v_{\mathrm{h},2}(\xi)=\frac{1}{12} \left(-6 \xi ^2 \cos (\xi )-12 \xi  \sin (\xi )+2 \pi
   ^2 \cos (\xi )-3 \cos (\xi )\right).
\]
Finally, substituting into the nonlinearity, evaluating at order $\nu^3$, we can directly solve for $\mu_3$ and find $
\mu_3=\frac{2\pi}{3}.
$ Since then $\elll=\pi-\pi\nu+\rmO(\nu^2)$ and $\mu=\frac{2\pi}{3}\nu^3+\rmO(\nu^4)$, we find the desired expansion \eqref{e:Luniv}.

Remarkably, all the coefficients computed here only depend on the first harmonic in the potential so that the expansion is indeed universal. 

%


\subsection{Vacuum formation with continuous repulsive potentials}\label{e:smooth}
We now turn to the case of a general repulsive potential. Unfortunately, a simple cosine ansatz is not possible in this case since the repulsive potential is not rank-one, that is, it does not map in the space spanned by constants and cosines. We therefore return to the free-boundary formulation involving a rescaling of potentials. Absorbing a smooth part of the potential into the attractive potential, we may assume without loss of generality that the kernel $V_\mathrm{rep}$ is given by the Green's function kernel to  $(1-\eta^2\partial_{xx})^{-\beta}$. Our free-boundary formulation \eqref{e:Fv}--\eqref{e:Fm} then becomes
\begin{align}
 0=F_v(v,\elll,\rho,\mu)[z]&:=\int_{-\pi}^\pi V^\elll_\mathrm{rep}(z-\xi)v(\xi)\rmd\xi -\int_{-\pi}^\pi V_\mu^\elll(z-\xi)v(\xi)\rmd\xi -\rho,\label{e:Fv1}\\
 0=F_\mathrm{bc}(v)&:=v(\pi),\label{e:Fbc1}\\
 0=F_m(v,\elll)&:=\frac{\elll}{2\pi^2}\int_{-\pi}^\pi v(\xi)\rmd\xi -1.\label{e:Fm1}
\end{align}
The rescaled kernel satisfies, in the case $\beta=1$, 
\[
\left(1-\frac{\eta^2\elll^2}{\pi^2}\partial_{xx}\right) V^\elll_\mathrm{rep}=\delta_0+\alpha(\elll,\eta)\delta_\pi,
\]
where $\alpha(\pi,\eta)=0$ stems from the discontinuity from truncating the potential at $L<\pi$ due to scaling. For $\beta<1$, the correction to the Dirac-$\delta$ at the origin is more regular, a distribution with Fourier coefficients $\sim |k|^{2(\beta-1)}$. We may therefore apply $
\left(1-\frac{\eta^2\elll^2}{\pi^2}\partial_{xx}\right) ^\beta$ to \eqref{e:Fv1} and obtain 
\begin{equation}\label{e:Fv2}
 0=\tilde{F}_v(v,\elll,\rho,\mu)[z]:=v(z) + (D(\eta,L,\beta)v)[z]-\int_{-\pi}^\pi \tilde{V}_\mu^\elll(z-\xi)v(\xi)\rmd\xi -\rho,
\end{equation}
with $D(\eta,L,\beta)$ continuous in $\eta$ and $L$ as an operator on $C^0$. The modified potential $\tilde{V}$ is continuous with our assumptions on smoothness of $V_\mathrm{att}$. From hereon, we simply follow the reasoning in \S\ref{s:fbdy}, including the parameter $\eta$.

\begin{figure}
\hfill\includegraphics[height=2in]{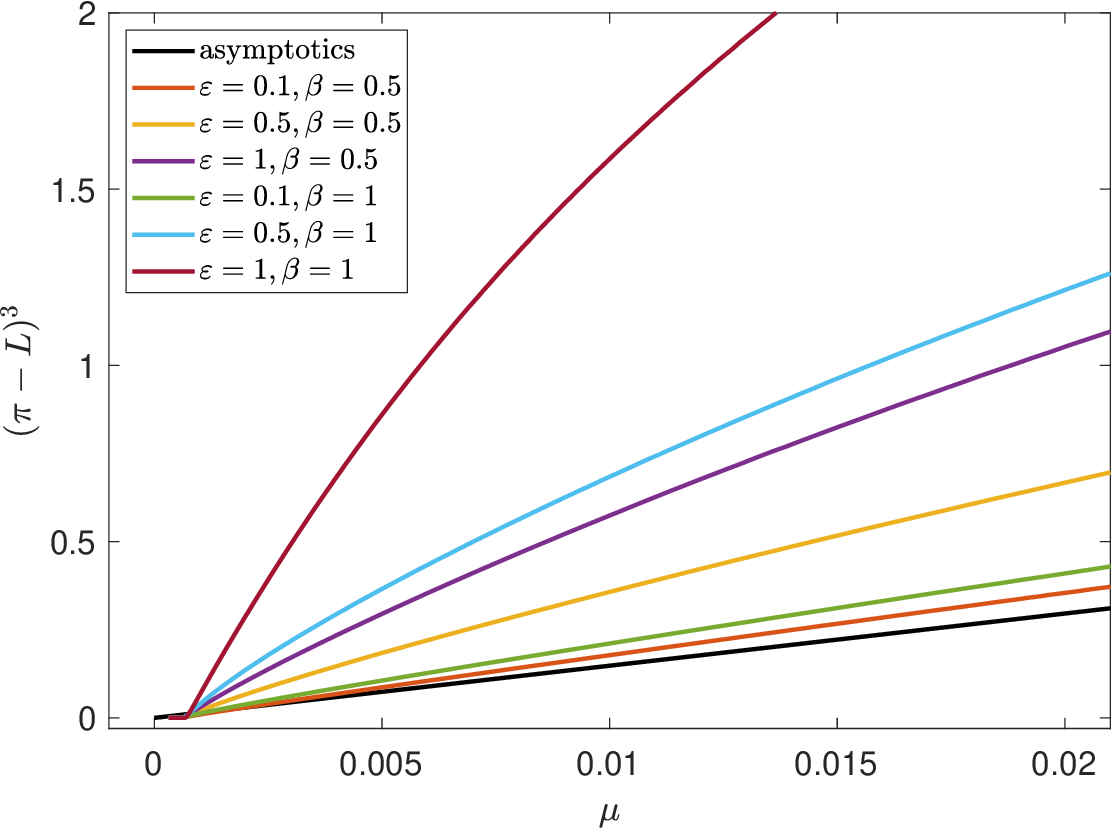}\hfill
\includegraphics[height=2in]{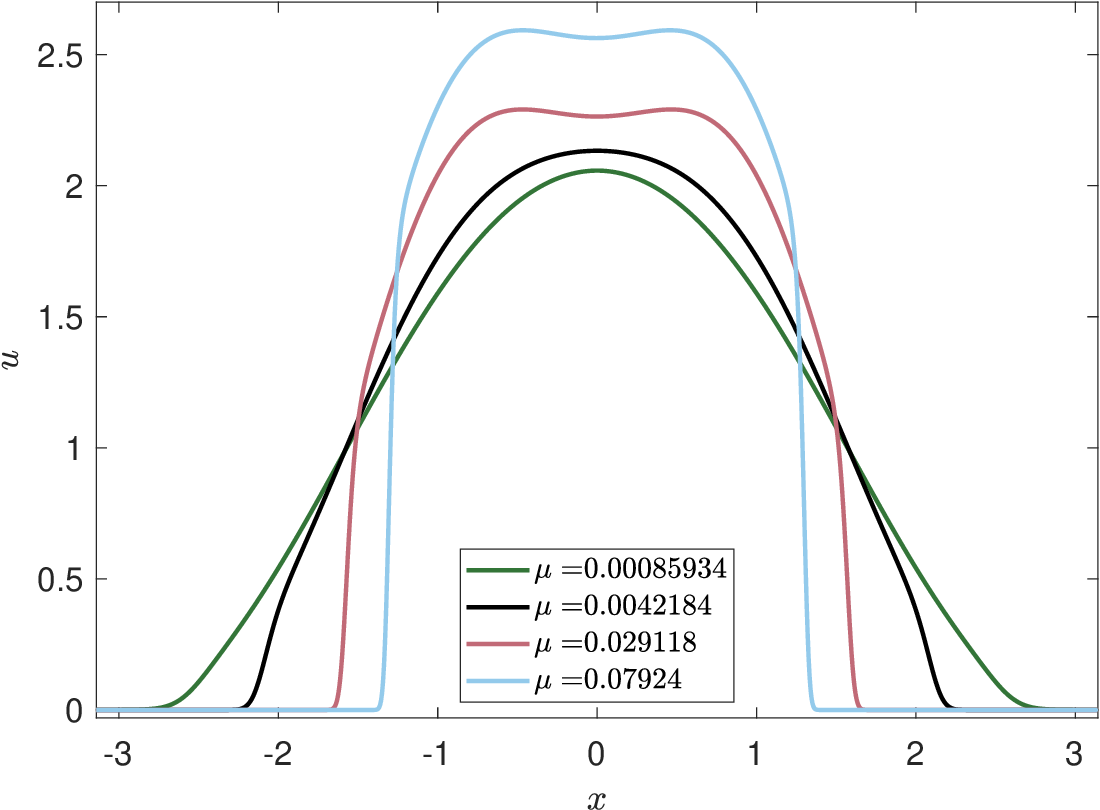}\hfill $ $
\caption{Computational results for vacuum formation with smoothed repulsive potential $(1-\eta^2\partial_{xx})^{-\beta}$ with small regularizing viscosity $10^{-3}$. The size of vacuum bubbles still scales with $\mu^{1/3}$ but coefficients depend on $\eta$ and $\beta$ (left); asymptotics for the case of Dirac-potential are also shown for reference. The shape of clusters changes significantly, with much steepened drop-off of concentrations near the boundary of the cluster and, for larger $\mu$, a dip at the center of the cluster (right).
}
\label{f:frac}
\end{figure}

The expansion for $\mu$ still vanishes at second order, but the third-order coefficient is no longer universal and depends on both $\eta$ and $\beta$. Fig. \ref{f:frac} shows the resulting change in the asymptotics and the formation of different cluster shapes as $\mu$ increases. In fact, the cubic expansion holds for all values of $\eta$, but we were not able to determine the direction of branching, in general. As a general rule, we confirm that weaker, that is, smoother repulsion modeled by increasing $\beta$ and/or $\eta$ gives larger vacuum regions as it allows for easier clustering. 

\section{Vacuum bubbles --- systems}
For systems, we write $u=(u_1,u_h)$ with $u_h\in \R^{P-1}$ collecting all components that do not form vacuum near $\mu=0$, and find the steady-state equations
\begin{align*}
u_1\cdot \left(V_{11}*u_1+V_{1h}*u_h \right)_x&=0,\\
\left(V_{h1}*u_1+V_{hh}*u_h \right)_x&=0,
\end{align*}
where we omitted the factor $u_h$ in the second equation since $u_h$ is nonzero. As a consequence, we can solve the second equation by inverting convolution with $V_{hh}$, which is possible since the Fourier transform of $V_{hh}$ is invertible. We write $V_{hh}^{-1}$ for the kernel of the inverse.  Substituting the result in the first equation, we find
\[
u_1\cdot \left(\tilde{V}_{11}*u_1 \right)_x=0,\qquad \tilde{V}_{11}=V_{11}+V_{1h}*V_{hh}^{-1}*V_{h1}.
\]
The resulting equation is thus precisely the same equation as in the one-dimensional case. When the repulsion is given by Dirac-$\delta$ potentials, we can again compute the expansion of the critical branch explicitly.

In the case when the principal part of $V(x)=a \delta(x) + b \cos(x)$ for $P\times P$-matrices $a$ and $b$, one can under certain conditions again find profiles of the form $(u^0 + u^1 \cos(x))_+$ with vectors $u^0=(u^0_{1},u^0_h)$ and $u^1=(u^1_{1},u^1_h)$. {For} consistency, we need to require that $a$ is block-diagonal, that is, $a_{j1}=a_{1j}=0$ for $j>1$. Writing the matrices in block form with blocks $a_{11},a_{1h},a_{h1}, a_{hh}$, and analogously for $b$, we find the reduced equation
\[
u^1_1=a_{11}^{-1}\left(b_{11}+\pi b_{1h}(a_{hh}-\pi b_{hh})^{-1}b_{h1}
\right) \int_{-L}^L \cos y (u^0_1+u^1_1 \cos y)\rmd y-\rho,
\]
which is analogous to the first equation in \eqref{e:vac1}. Adding the mass constraint and boundary condition at $x=L$, one can easily solve the resulting equations for $u^1_1$ and $u^0_1$. As expected, the agreement with the asymptotics is very similar to the scalar case shown in Fig. \ref{f:lmu}. We emphasize however that this rank-one approximation breaks down quickly as other species develop vacua as well.

We also studied \eqref{e:sys} numerically, using numerical continuation with artificial viscosity $\eps=0.03$.
We chose strong short-range repulsion between all particles, $a_{ij}>0$ for $i,j\in\{1,2\}$ and long-range self-repulsion $b_{11},b_{22}<0$. We then chose the strength and sign of the inter-species long-range interaction as a parameter $\kappa$, $b_{12}=b_{21}=\kappa$ so that $\kappa>0$ corresponds to mutual attraction and $\kappa<0$ to mutual repulsion. Specifically, we let
\begin{equation}\label{e:coeff}
\left(\begin{array}{cc}
a_{11}& a_{12}\\
a_{21}& a_{22}
\end{array}\right)
=
\left(\begin{array}{cc}
0.8& 1\\
1& 1
\end{array}\right),\qquad
\left(\begin{array}{cc}
b_{11}& b_{12}\\
b_{21}& b_{22}
\end{array}\right)
=
\left(\begin{array}{cc}
-0.3& \kappa\\
\kappa& -0.3
\end{array}\right).
\end{equation}
One finds the bifurcation points $\kappa_\mathrm{jc}=0.8403$ where joint clustering sets in for $\kappa>\kappa_\mathrm{jc}$ and $\kappa_\mathrm{seg}=-0.3310$ where segregation sets in for $\kappa<\kappa_\mathrm{seg}$. Note that segregation is facilitated with a much weaker mutual interaction strength than clustering as it needs to overcome significantly less short- and long-range repulsion.

\begin{figure}
\hfill\includegraphics[height=2in]{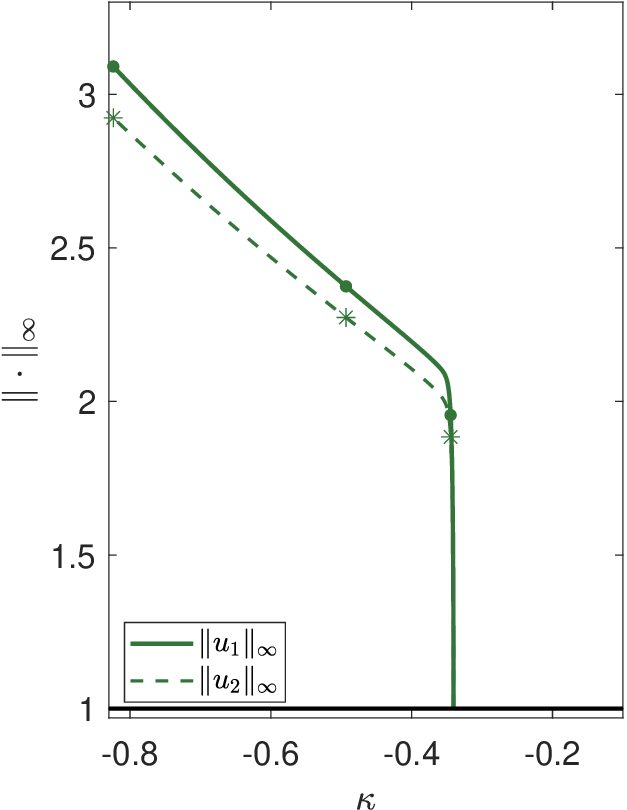}\hfill
\includegraphics[height=2in]{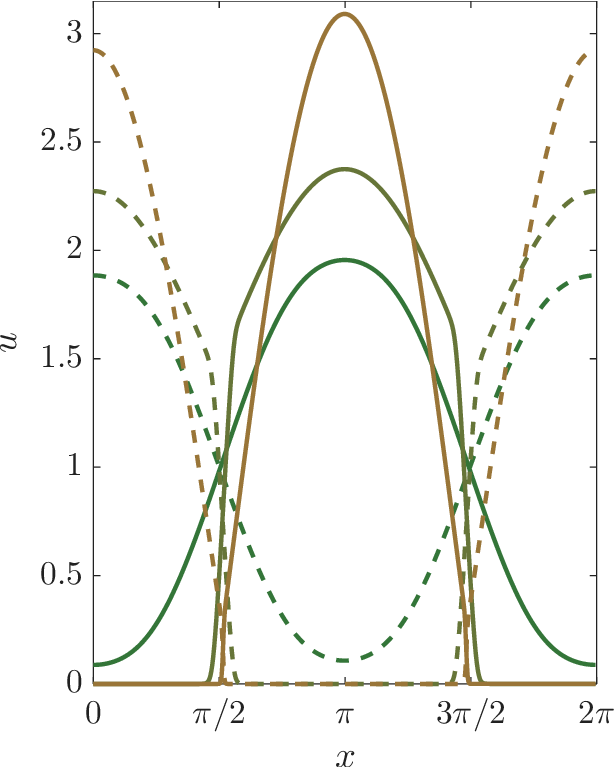}
\hfill\includegraphics[height=2in]{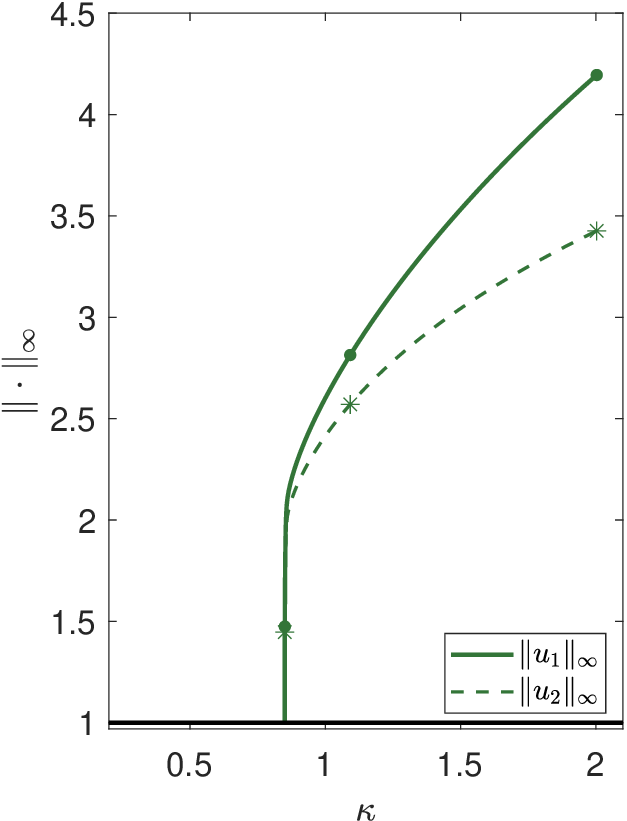}\hfill
\includegraphics[height=2in]{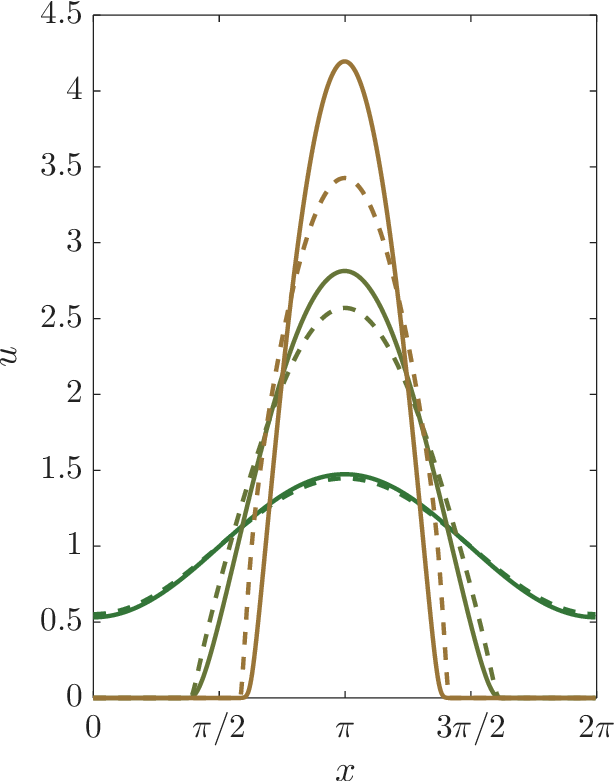}\hfill
\caption{Numerically computed bifurcation diagrams for \eqref{e:sys} with coefficients \eqref{e:coeff}. For $\kappa<0$ decreasing, the inter-species repulsion leads to segregation (left two panels), while increasing $\kappa>0$ leads to joint clustering (right two panels). First species is always shown as solid lines, second as dashed lines. Sample profiles are shown at the marked points with increasing amplitude as the branch is followed. Throughout, we added diffusion with strength $\eps=0.03$.
}
\label{f:sys}
\end{figure}

We confirmed the asymptotics for the size of vacuum regions from Thm.~\ref{t:1}.
Bifurcation diagrams and sample profiles are shown in Fig. \ref{f:sys}.

\section{The effect of diffusion --- a scalar case study}
\begin{figure}
\hfill\includegraphics[height=2in]{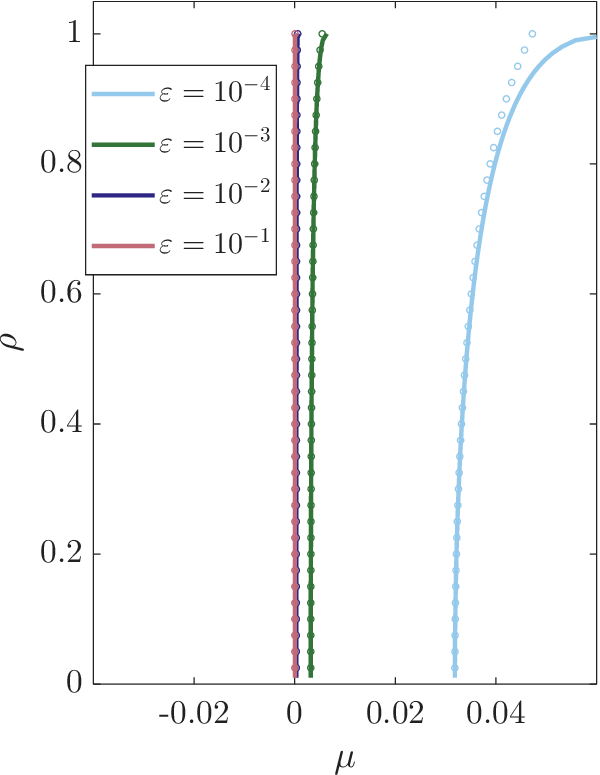}\hfill
\includegraphics[height=2in]{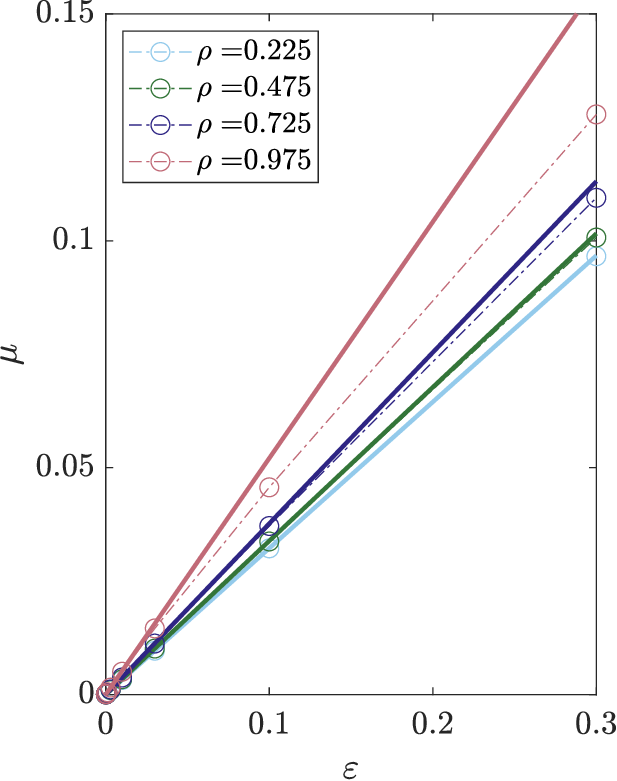}
\hfill\includegraphics[height=2in]{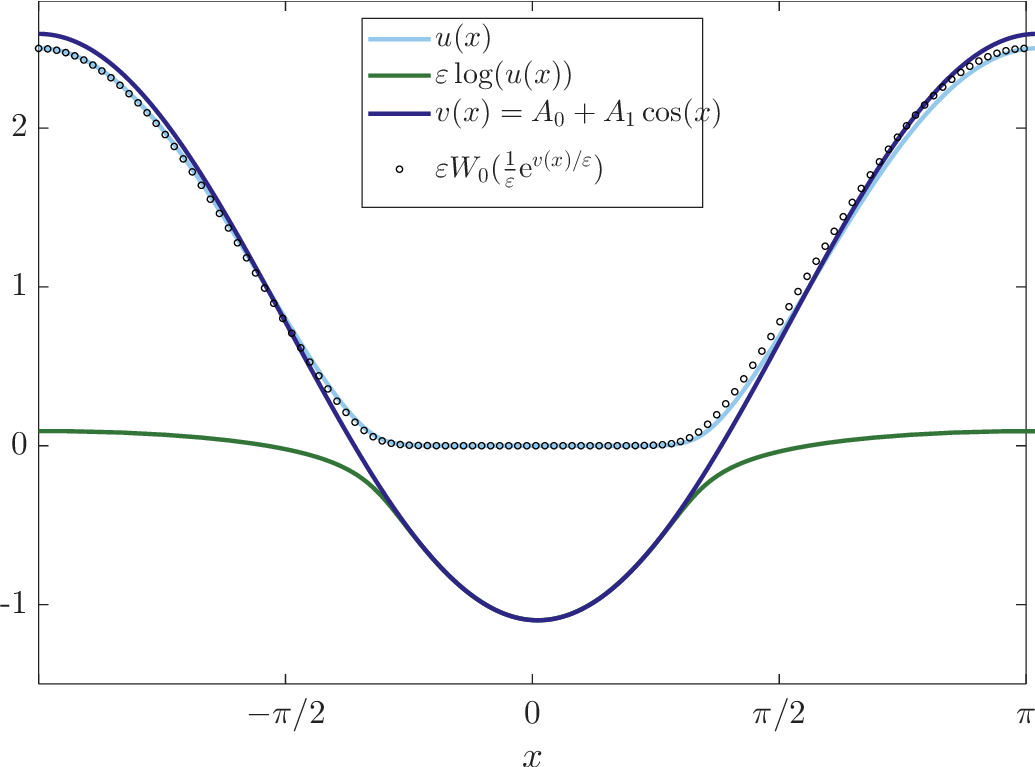}\hfill
\caption{Corrections from diffusion to bifurcation diagram and solutions to \eqref{e:vv}: the vertical branch bends with an $\rmO(\eps)$-correction as in \eqref{e:epsmu}, $\rho$ as in \eqref{e:norm} (left, computed circles, predicted solid line); shift $\mu$ of the branch at fixed $\rho$ for varying $\eps$  (center, computed circles, predicted solid line); profiles with amplitude $\|u\|_\infty=2.5$ at $\eps=0.1$, compared with the prediction from \eqref{e:uexpl} and \eqref{e:vexpl} and asymptotics \eqref{e:asyw} (right; fitted $A_1$ from curvature at minimum).}
\label{f:eps}
\end{figure}

A natural question  in light of applications such as the lipid rafts discussed in the {introduction} is the robustness of the non-hysteretic switching when the system is subjected to noise. Adding small noise to particle dynamics, one finds in a continuum limit an additional diffusion term,
\begin{equation}\label{e:vv}
u_t=\eps u_{xx}+(u\cdot(u-(\frac{1}{\pi}+\mu)V*u)_x)_x,\qquad V(x)=\cos(x),
\end{equation}
with $\eps\gtrsim 0$. At $\mu=0$ and $\eps=0$, \eqref{e:vv} has a vertical branch of steady-states $u_*=1+\rho\cos(x)$. The next result is that noise induces an $\rmO(\eps)$-hysteresis, that is, the vertical branch is weakly slanted into a supercritical pitchfork bifurcation.
\begin{proposition}[Diffusive  corrections]
The diffusive Vlasov equation \eqref{e:vv} possesses an almost vertical, $2\pi$-periodic, even branch of solutions
$u_*(x;\rho,\eps)$, $0<\rho<1$, for $\mu=\mu_*(\rho,\eps)$ with expansion
\begin{equation}
\mu_*(\rho,\eps)=\eps\mu_1(\rho)+\rmO(\eps^2),\quad u_*(x;\rho,\eps)=1+\rho\cos(x)+\rmO(\eps),
\end{equation}
and normalization
\begin{equation}\label{e:norm}
\int_{-\pi}^\pi u_*(x)\rmd x=2\pi,\qquad \int_{-\pi}^\pi \cos(x)u_*(x;\rho,\eps)\rmd x=\rho\pi.
\end{equation}
Explicitly, we have
\begin{equation}\mu_1(\rho)=\frac{2}{\pi\rho^2}\left(1-\sqrt{1-\rho^2}\right),
\end{equation}
with limits $\mu_1(0)=1/\pi$ and $\mu_1(1)=2/\pi$.
\end{proposition}
The results  predict very accurately the numerically computed almost vertical part of the bifurcation branches in Fig. \ref{f:2a}, left panel. Detailed comparisons are shown in Fig. \ref{f:eps}.

One can of course analyze the pitchfork bifurcation for fixed $\eps$ and then finds, in addition to the shift of the bifurcation point to $\mu_*=\eps$, an $\rmO(\eps)$ cubic {coefficient} in the pitchfork bifurcation of the form $\eps\mu_1''(0)/2$. For fixed $\eps>0$, bifurcations of this type have been studied in the literature; see for instance \cite{CGPS}. We emphasize that the bifurcation results obtained in this fashion only allow for $|\rho|\ll 1$, however, thus do not capture the effectively more dramatic almost-absence of hysteresis when $\eps\ll 1$.

\begin{Proof}
We write the steady-state equation abstractly as
\begin{equation}\label{e:stsys}
\begin{aligned}
F(u,\mu,\eps)&:=\eps u_{xx}+u\cdot(u-(\frac{1}{\pi}+\mu)V*u)_x)_x=0,\\
F_0(u)&:=\int_{-\pi}^\pi u(x)\rmd x-2\pi,\\
F_1(u)&:=\int_{-\pi}^\pi \cos(x)u(x)\rmd x-\rho\pi.
\end{aligned}
\end{equation}
Altogether, this defines a map
\[
G(u,\mu,\eps):H^2_\mathrm{e,p}\times \R^2\to \mathring{L}_\mathrm{e,p}^2\times \R^2,
\]
where the subscripts $\{\mathrm{e,p}\}$ recall the fact that we are working in even, periodic functions, and $\mathring{L}^2$ denotes functions with zero average. Here and in the following we think of $\rho$ fixed with $|\rho|<1$. One easily verifies that $G$ is well-defined and smooth, and
\[
G(u_*^0(\cdot;\rho),0,0,\rho)=0,\qquad u_*^0(x;\rho)=1+\rho \cos(x).
\]
The derivative at this zero is given through
\[
\mathcal{A}=\left(\begin{array}{ccc}
\mathcal{L} & \partial_\mu F_* &  \partial_\eps F_*\\
\langle 1,\cdot\rangle & 0& 0\\
\langle \cos(\cdot),\cdot\rangle & 0 & 0
\end{array}\right),
\]
where
\begin{equation}\label{e:terms}
\begin{aligned}
\mathcal{L} u&=(u_*^0\cdot(u-\frac{1}{\pi}V*u)_x)_x,\\
\partial_\eps F_*&=-\rho \cos(\cdot),\\
\partial_\mu F_*&=\pi\rho\cos(\cdot)+\rho^2\pi \cos(2x),
\end{aligned}
\end{equation}
and $\langle 1,\cdot\rangle$ and $\langle \cos(\cdot),\cdot\rangle$ are continuous linear functionals on $H^2_\mathrm{e,p}$.

The operator $\mathcal{L}$ is an strictly elliptic operator due to the fact that we restrict to $|\rho|<1$, and therefore Fredholm with index 1 due to the restriction of the codomain to average zero functions. Bordering lemmas for Fredholm operators then give that $\mathcal{A}$ is Fredholm of index 1, so that we expect a one-dimensional set of solutions, which turns out to be the branch described in the result with $\rho$ fixed.
We next move on to identify kernel and cokernel of $\mathcal{A}$.

The kernel of $\mathcal{L}$ is spanned by $\{1,\cos(\cdot)\}$, thus implying a one-dimensional cokernel. In order to identify the cokernel explicitly, note that
\[
\mathcal{L}=\mathcal{M}\mathcal{L}_0,\qquad \mathcal{M}v=\partial{x}(u_*^0\partial_xv),\qquad \mathcal{L}_0 u=u-\frac{1}{\pi} V*u,
\]
so that the adjoint is $\mathcal{L}^*=\mathcal{M}\mathcal{L}_0$ due to the fact that both $\mathcal{M}$ and $\mathcal{L}_0$ are self-adjoint. Solving $\mathcal{M}v=\cos(\cdot)$, we find the kernel of the adjoint
\[
e_0^*(x)=-\frac{1}{\rho} \log(1+\rho\cos(x)).
\]
As a consequence, we find that $\partial_\mu F_*\not\in\rg(\mathcal{L})$, through
\begin{equation}\label{e:muproj}
\langle e_0^*,\partial_\mu F_*\rangle = \int_{-\pi}^\pi \frac{-1}{\rho} \log(1+\rho\cos(x)) \left(\pi\rho\cos(\cdot)+\rho^2\pi \cos(2x)\right)=-\pi^2\rho\neq 0.
\end{equation}
Inspecting the first two columns of $\mathcal{A}$, we find that $\mathcal{A}(u,\mu,0)=0$ if and only if $u=0$ and $\mu=0$, as the second and third component of $\mathcal{A}$ imply $u\not\in\ker\mathcal{L}$ and therefore $\mathcal{L}u\neq 0$, $\mathcal{L}u\neq \partial_\mu F_*\mu$ unless $u=0$ and $\mu=0$. We can therefore solve for $u,\mu$ as functions of $\eps$ and obtain the leading order expansion of $\mu$ by projecting the first equation onto the cokernel. Evaluating also
\begin{equation}\label{e:epsproj}
\langle e_0^*,\partial_\eps F_*\rangle = \int_{-\pi}^\pi \frac{-1}{\rho} \log(1+\rho\cos(x)) \left(-\rho\cos(x)\right)=\frac{2\pi}{\rho}\left(1-\sqrt{1-\rho^2}\right)\neq 0,
\end{equation}
this gives
\begin{equation}\label{e:epsmu}
\mu=\mu_1 \eps+\rmO(\eps^2),
\qquad \mu_1=-\frac{ \langle e_0^*,\partial_\eps F_*\rangle   }{\langle e_0^*,\partial_\mu F_*\rangle}=\frac{2}{\pi\rho^2}\left(1-\sqrt{1-\rho^2}\right).
\end{equation}
\end{Proof}
Clearly, the analysis so far is not reliant on the specific form of the potential and could also be carried out for more general repulsion kernels and systems. In the remainder of this section, we exploit more explicitly the fact that the potential is rank-one to extract asymptotics in the case where the solution possesses a vacuum region.

We can write the stationary form of \eqref{e:vv} after integrating as
\[
\eps u_x + u\cdot (u-(\frac{1}{\pi}+\mu)A \cos(x))_x=0,\qquad A=\int_{-\pi}^\pi \cos(x)u(x)\rmd x.
\]
The first equation can be explicitly solved setting $v=\eps \log u + u$, or
\begin{equation}\label{e:uexpl}
u=\eps W_0(\frac{1}{\eps}\rme^{v/\eps}),
\end{equation}
where $W_0$ is the first, positive branch of the $W$-Lambert function. The resulting differential equation for $v$ then gives
\begin{equation}\label{e:vexpl}
v(x)=-A(\frac{1}{\pi}+\mu)\cos(x)+m.
\end{equation}
Exploiting asymptotics of the $W$-Lambert function, or simply asymptotics on the equation for $v$ in terms of $u$, we find that
\begin{equation}\label{e:asyw}
u(x)=\left\{\begin{array}{ll}
v(x)+\rmO(\eps), \quad & \text{ when }v(x)>0,\\
\rme^{-|v(x)|/\eps}(1+\rmO(\rme^{-|v(x)|/\eps}/\eps).
\end{array}\right.
\end{equation}
As a consequence, profiles of the form $(A_0+A_1\cos(x))_+$ inherit $\rmO(\eps)$ corrections in the regime where they are positive, and they are exponentially small in $\eps$ away from the boundary of the vacuum region.

Summarizing, in this scenario we have an almost explicit shape given through the W-Lambert function,
\eqref{e:uexpl} and \eqref{e:vexpl}, with expansions for the profile away from the boundary of the vacuum \eqref{e:asyw}. The results give excellent predictions for numerical {computations} as demonstrated in Fig. \ref{f:eps}.

\section{Discussion}\label{s:d}
We discussed a curious vertical bifurcation in limits of many-particle systems. Key to the phenomenon is a destabilization of a crystalline, uniform state due to a long-range attractive force that outcompetes a strong short-range attraction at intermediate distances, leading to the formation of clusters and gaps, and also to segregation when several particle species interact with each other. Quite illuminating is the study of a special rank-one case of periodized attractive potentials being simply a multiple of a cosine. We also discussed the effects of noise by adding a diffusive term to the equation. Technically, our main result describes the opening of a vacuum bubble with a universal scaling in the parameter.

There are clearly numerous questions that arise from this work.

\paragraph{Dynamics. } We only describe the bifurcation diagram, here. It would be interesting to describe additional features, starting with linear stability, and hopefully covering, more interestingly, the existence of slow manifolds in the sense of Fenichel's geometric singular perturbation theory \cite{fenichel} for $\mu\sim 0$. Our analysis of robustness under addition of small diffusion exploited the fact that the vertical branch has an invertible linearization, up to the tangent vectors induced by the vertical branch, translations, and changes in total mass. Given the gradient-flow structure of the equation, we suspect that this invertibility is equivalent to normal hyperbolicity of the vertical branch of equilibria  in the sense of \cite{fenichel}. Before vacuum formation,  this hints at the presence of an asymptotically attracting slow manifold that mediates the slow formation of vacuum for $\mu\gtrsim 0$ and the collapse of vacuum bubbles for $\mu\lesssim 0$. The main technical difficulty here appears to overcome the lack of normal hyperbolicity at the moment of vacuum formation. Clearly, this result would more conclusively establish the absence of hysteresis that we infer here only statically.

\paragraph{Finite-size effects.} The Vlasov-limit is usually obtained only in a weak sense as $N\to\infty$. Nevertheless, one may hope for a more detailed description, at least in one space-dimension. We performed numerical continuation in finite-particle interaction with the repulsive potential from \eqref{e:gauss} and a cosine attractive potential, finding bifurcation diagrams with several branches, which nevertheless display striking closeness to the continuum limit.
\begin{figure}
\includegraphics[width=0.49\textwidth]{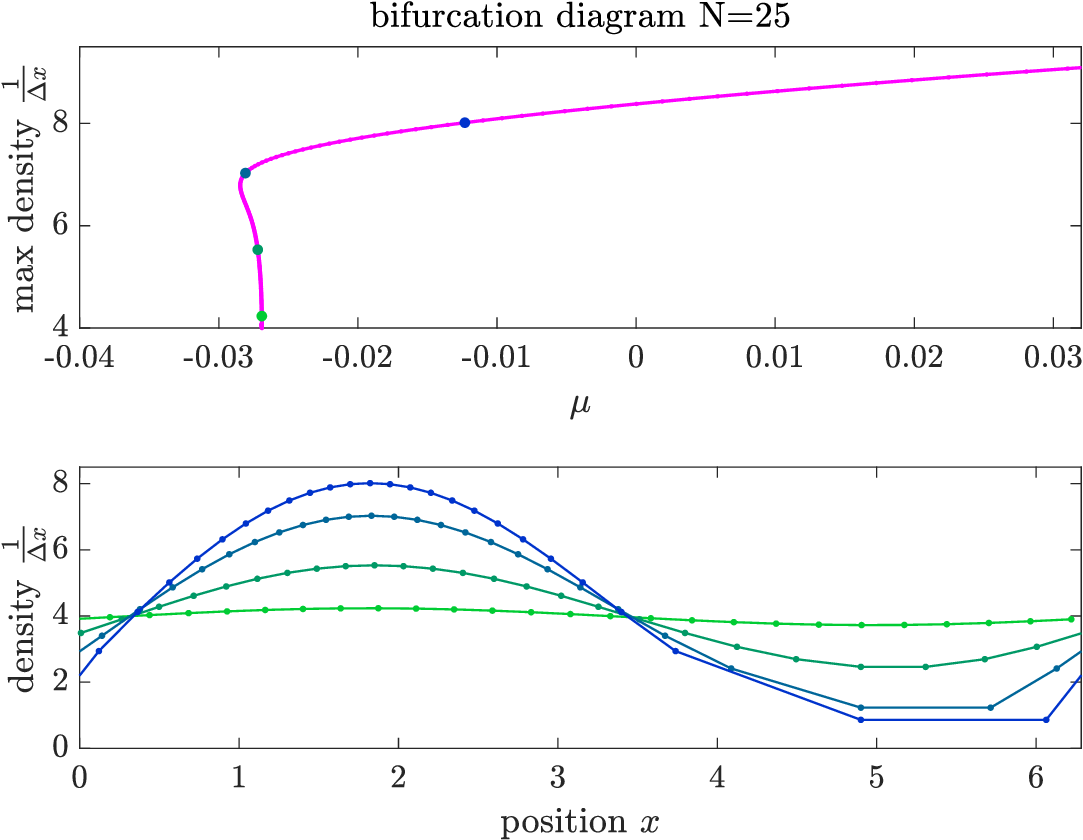}\hfill
\includegraphics[width=0.49\textwidth]{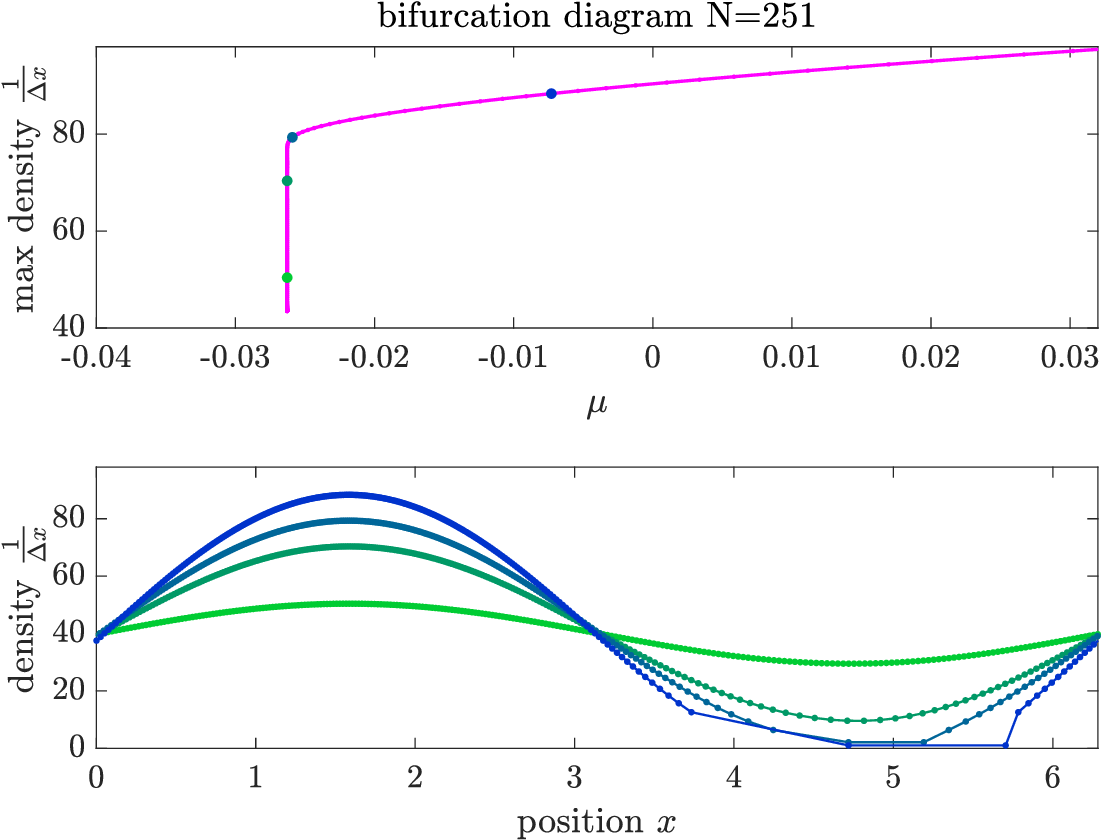}
\caption{Bifurcation diagram of maximal density versus bifurcation parameter in case of finite particle sizes $N=25,251$ (top left and top right) in equation \eqref{e:partper} with repulsive potential from \eqref{e:gauss}, $\delta=0.3$, and an attractive cosine potential $-(\frac{1}{\pi}+\mu)\cos(x)$. Bottom panels show position of particles at marked parameter values plotted together with the inverse distance to the neighboring particle as a proxy for the local density.}
\label{f:2}
\end{figure}

\paragraph{Higher space-dimension.} We study higher-dimensional particle interaction in a forthcoming paper \cite{swarm2d_reu}. In addition to essentially one-dimensional vacuum regions, given as stripes in a periodic grid, we then find roughly spherical bubbles, also with universal expansions for the size in terms of the parameter, with an interesting competition between vacuum bubbles and stripes. One would expect finite-size effects to be significantly more complex to the non-monotonicity of
rearrangements in cluster and vacuum formation.

\paragraph{Reversible switching elsewhere.} We noted in the introduction that the absence of hysteresis is rather unusual when starting from generic low-dimensional bifurcations. It does not appear to arise near generic low-codimension bifurcation points, even in the presence of symmetry. Bifurcation branches similar to the finite-size diagrams shown in Fig. \ref{f:2} arise in singular Hopf bifurcations, known there as canard explosions; see for instance \cite{KScanard,BE86}, and references therein.

In a quite different context, almost reversible bifurcations arise with generic instabilities in large domains: the instability is mediated by a front that moves into the domain for $\mu>0$ and out of the domain for $\mu<0$. Weak interaction with the boundary then leads to almost vertical bifurcation diagrams; see \cite{ADSS}.

We emphasize that there is no apparent link with these different settings other than that they all permit an almost reversible switching.

\paragraph{Lipid rafts revisited.} Circling back to our motivation by biological switching, we identified here possible mechanisms and
suitable mathematical models that can describe
fast and robust switching, back and forth between different functional patterns, particularly in biological settings. The example we discussed in the introduction are so-called lipid rafts.
Individual lipid molecules are able to diffuse
two-dimensionally within the cellular lipid bilayer, and in artificial lipid bilayers 
one can observe domains of different lipid composition.
Such specialized domains are thought to be involved in signal transduction 
across the cell membrane and in organizing membrane trafficking, as well as in the organization 
and regulation of membrane proteins. One major driver for this
segregation phenomenon could be the selective association between certain 
lipids. Our analysis isolates this selective association in the study of a continuum limit for interacting particle systems, driven solely by attraction, repulsion, and possibly diffusion. We do then indeed demonstrate that near parameter values of balance between inner- and intra-species attractive and repulsive forces, made precise in an explicit linear analysis near a perfectly mixed state, reversible switching is possible: small changes in parameter values, for instance sensitivity of particles to interaction forces, can drive the system from perfectly mixed states to almost perfectly sorted states or states exhibiting vacua --- and reverting the small parameter change will restore the perfectly mixed state.
\bibliographystyle{abbrv}

\end{document}